\documentclass[reqno]{amsart}
\usepackage{amsfonts,amssymb}
\usepackage{graphicx} 
\usepackage{float}
\usepackage{color}

\usepackage[rightcaption]{sidecap}
\usepackage{caption} 
\usepackage{soul}
\usepackage{dsfont}
\usepackage{mathrsfs}

\newcommand{\stkout}[1]{\ifmmode\text{\sout{\ensuremath{#1}}}\else\sout{#1}\fi}

\usepackage{amsthm, amssymb, latexsym, amsmath, bm 
}
\usepackage{ulem} 
\usepackage[margin=1.35in]{geometry}
\usepackage{enumerate}
\usepackage[toc,page]{appendix}

\usepackage{dsfont}

\usepackage{hyperref}
\usepackage{cleveref}
\usepackage{xcolor}
\usepackage{subcaption}

\usepackage{tikz}

\hypersetup{
  colorlinks   = true, 
  urlcolor     = blue, 
  linkcolor    = blue, 
  citecolor   = red 
}

\newtheorem{theorem}{Theorem}[section]

\newtheorem{corollary}[theorem]{Corollary}

\newtheorem{lemma}[theorem]{Lemma}

\newtheorem{remarks}{Remarks}
\newtheorem{remark}[remarks]{Remark}
\newtheorem{proposition}[theorem]{Proposition}
\theoremstyle{definition}

\numberwithin{theorem}{section}
\numberwithin{equation}{section}

\newcommand{\nc}{\newcommand}

\nc{\be}{\begin{equation}}
\nc{\la}{\label}
\nc{\ba}{\begin{array}}
\nc{\ea}{\end{array}}
\nc{\bs}{\begin{split}}
\nc{\es}{\end{split}}

\nc{\eps}{\epsilon}
\nc{\e}{\epsilon}
\nc{\lam}{\lambda}
\nc{\G}{\Gamma}
\nc{\g}{\gamma}
\nc{\al}{\alpha}
\nc{\del}{\delta}
\nc{\Del}{\Delta}\nc{\Om}{\Omega}
\nc{\Omt}{\tilde{\Omega}}
\nc{\ta}{\tau}
\nc{\w}{\omega}
\nc{\io}{\iota}
\nc{\z}{\zeta}
\nc{\s}{\sigma}
\nc{\Si}{\Sigma}
\nc{\Lam}{\Lambda}

\newcommand{\R}{{\mathbb R}}

\newcommand{\bS}{{\mathbb S}}

\newcommand{\cS}{\mathcal{S}}

\nc{\ra}{\rightarrow}

\nc{\ran}{\rangle}
\nc{\lan}{\langle}

\newcommand{\one}{{\bf 1}}

\newcommand{\n}{\nabla}
\newcommand{\p}{\partial}

\newcommand{\DETAILS}[1]{}

\newcommand{\variants}[1]{}

\newcommand{\bu}{\textbf{u}}

\newcommand{\Dr}{\Delta^{\!\scriptscriptstyle\operatorname{rad}}}

\newcommand{\Ar}{A^{\!\scriptscriptstyle\operatorname{rad}}}

\newcommand{\mcV}{{\mathcal V}}
\newcommand{\mcF}{{\mathcal F}}
\newcommand{\mcN}{{\mathcal N}}
\newcommand{\mcM}{{\mathcal M}}
\newcommand{\mcR}{{\mathcal R}}

\begin{document}
\null\hfill\begin{tabular}[t]{l@{}}
\scriptsize\textbf{Journal reference:} Journal of Nonlinear Science \\
\scriptsize \textbf{DOI:} 10.1007/s00332-025-10152-9\\ \vspace{0.1cm}
\end{tabular}

\vspace{-0.15cm}

\title[]{The FitzHugh-Nagumo System on Undulated Cylinders: Spontaneous Symmetrization and Effective System}

\author[]{G. Karali$^{\dag,\#}$,\, K. Tzirakis$^{\ddag,\#}$,\, I. M. Sigal$^{*}$}

\thanks{$^{\dag}$ Department of Mathematics, National and Kapodistrian University of Athens, Panepistimiopolis 15784, Athens, Greece.}
\thanks{$^{\#}$ Institute of Applied and Computational Mathematics, FORTH, 71110, Heraklion, Greece.}
\thanks{$^{\ddag}$ Department of Mathematics and Applied Mathematics, University of Crete, 71409, Heraklion, Greece.}
\thanks{$^{*}$ Department of Mathematics, University of Toronto, Toronto, ON M5S 2E4, Canada.}


\begin{abstract}
We consider the FitzHugh-Nagumo system on undulated cylindrical surfaces modeling nerve axons. We show that for sufficiently small radii and for initial conditions close to radially symmetrical ones, (i) the solutions converge to their radial averages, and (ii) the latter averages can be approximated by solutions of a $1+1$ dimensional (`radial') system (the \textit{effective} system) involving the surface radius function in its coefficients. This perhaps explains why solutions of the original $1+1$ dimensional FitzHugh-Nagumo system agree so well with experimental data on electrical impulse propagation.
\end{abstract}



\maketitle

\section{Introduction}\label{sec:intro}
The FitzHugh-Nagumo (FHN) system describes propagation of electrical impulses (pulses) along nerve filaments (axons). It is a simplification of the Nobel prize winning Hodgkin-Huxley (HH) system with slow unknown variables of ion concentrations, whose movement through the axon surface creates the potential gradient, clumped into a single slow variable. Despite this simplification, the FHN system provides a remarkably accurate qualitative and in many cases quantitative description of the rather complicated phenomenon of signal propagation along nerve cells. For comparison between the HH and FHN pulse profiles (as functions of a real variable $z$ representing a moving frame coordinate $x-ct$), see \cite{PhilSchuster} and Fig.~\ref{fig:comp} below.
\begin{SCfigure}[2][h]
\includegraphics[width=17em]{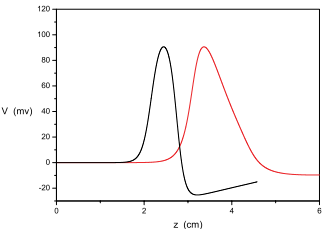}
\caption{Comparison of the FHN and HH pulse shapes (black and red curves, respectively) from \cite{PhilSchuster}. In both cases, the pulse speed is 1873 cm/sec and the pulse maximum is 90.6 mV (\cite{PhilSchuster}). (This figure is used with permission of World Scientific Publishing Company.)\\\phantom{k}}
\label{fig:comp}
\end{SCfigure}

Though axons are thin tubes of variable radii and of complicated geometry, with electrical impulses propagating on their surfaces, starting with the Hodgkin and Huxley work \cite{HodHux}, they are modeled by straight lines with no internal structure. To address this discrepancy, it was proposed by Talidou et al. \cite{TBS} to model axons by thin cylindrical surfaces and to extend the FitzHugh-Nagumo system (FHN) to such a surface, $\cS$, as
\refstepcounter{equation}\label{FHN-S1}
\begin{equation}\label{FHN-S}
\begin{split}
\frac{\partial u_1}{\partial t} &= \Delta_\cS u_1 + f(u_1) - u_2,\\  
\frac{\partial u_2}{\partial t} &= \eps(u_1-\gamma u_2),
\end{split}\tag{\theequation$\cS$}
\end{equation}

\noindent where $\Delta_\cS$ is the Laplace-Beltrami operator on a surface $\cS$. 
The function $f$ is the cubic polynomial 
\begin{equation}\label{nonlinf}
f(v) := -v (v - \al) (v - 1)\,,
\end{equation}
for $\al \in (0, \frac{1}{2}). $ The constants $\eps$ and $\gamma$ are positive and sufficiently
small, $0 <  \eps, \gamma \ll 1$. Apart from replacing $\partial^2_x $ by $\Delta_\cS, $ the FHN system is not changed. Thus, for $ \cS=\mathbb R, $ and consequently $ \Delta_\cS = \partial_x^2, $ Eq.~\eqref{FHN-S} is reduced to the standard FitzHugh-Nagumo system (\ref{FHN-S1}$\R$), in our notation, on $ \mathbb R $ (which does not involve $\theta$).

Specifically, \cite{TBS} considers \eqref{FHN-S} on cylindrical surfaces, $\cS$, called there the {\it warped cylinders}, whose radii vary slowly along their length and with the change of the azimuthal angle, defined as a graph over the standard cylinder, $\R\times \bS^1$, 
\begin{align}\la{graph-cyl}\cS_\rho:=\{(x, y, z) \vert r:=\sqrt{y^2+z^2} =\rho(x, \theta)\},\end{align}

\noindent where $\theta$ is the polar angle and $\rho(x, \theta)$ is a slowly varying function of $ (x, \theta) \in \R\times \mathbb S^1. $ The only difference with the standard FHN system (\ref{FHN-S1}$\R$) is in replacing the second derivative by the surface Laplacian.

Besides providing the next natural step in modeling pulse dynamics in axons, Eq.~\eqref{FHN-S} is an interesting mathematical model for propagation of traveling waves on surfaces.

In this paper, we study general properties of solutions for very thin cylindrical surfaces and emergence of 1D effective systems. To explain the issues involved, we consider Eq.~\eqref{FHN-S} first on $\cS = \cS_R, $ where $\cS_R$ is the surface of a straight cylinder, $\cS_R := \R\times R S^1$, of a constant radius $R$ centered about the $x$-axis in $\R^3$. We parametrize $ \cS_R $ by the cylindrical coordinates as
\begin{equation} \label{eq:def-SR}
\cS_R = \left \{ (x,R\cos\theta, R\sin\theta) 
\in\R^3\ \big\vert\ 
x\in \R,\, \theta\in [0,2\pi)
\right\}\, ,
\end{equation}
with the Riemannian area element and the Laplacian given by $ d\mu_R=R\,d\theta dx$ and
\begin{equation} \label{Delta-R}
\Delta_{\cS_R}= \partial_x^2 + R^{-2}\partial_\theta^2\,.
\end{equation}
Clearly, each solution of Eq.~(\ref{FHN-S1}$\mathbb R$) defines a radial, i.e.\ axi-symmetric, or $\theta$-independent, solution of Eq.~\eqref{FHN-S} for $\cS = \cS_{R},$ for any $R>0$ and vice versa.

Like Eq.~(\ref{FHN-S1}$\R$), Eq.~\eqref{FHN-S} for $\cS = \cS_{R},$ denoted (\ref{FHN-S1}$\cS_{R}$), is invariant under translations, i.e., if
$\bu(x, \theta, t)=(u_1(x, \theta, t), u_2(x, \theta, t))$ is a solution, then so are its translates
$$
\bu_h(x, \theta, t) := \bu(x - h, \theta, t)\,,\quad h\in \R.
$$
Hence Eq.~(\ref{FHN-S1}$\cS_R$) might have traveling wave solutions.
It is known that for $0<\eps\ll 1$, Eq.~(\ref{FHN-S1}$\R$) has two traveling wave solutions, fast and slow, the first of which is stable while the second unstable (see \cite{Flores,EvansIV,IMT87}). 

More precisely (see \cite{TBS}), for $0<\alpha<\frac12$ and $\eps$, $\gamma$ are positive and sufficiently small, Eq.~(\ref{FHN-S1}$\R$) has at least two different pulse solutions, the \textit{fast} pulse (shown in Fig.~\ref{fig:comp}), which travels with speed $c_f(\eps)= \frac{\sqrt{2}}{2}(1-2\alpha)+o(\eps)$, and a \textit{slow} pulse that travels with speed $c_s(\eps) = O(\sqrt{\eps})$ (Hastings~\cite{Hastings_76, Hastings_82}, Carpenter~\cite{Carpenter_77}, Langer~\cite{Langer_80}, Jones, Kopell and Langer~\cite{JKL91}, Krupa, Sandstede and Szmolyan~\cite{KSS97}, Arioli and Koch~\cite{AK15}).  The uniqueness of the fast pulse was proven by Langer~\cite{Langer_80}. Furthermore, the fast pulse is stable (Jones \cite{Jones_84} and  Yanagida~\cite{Yanagida_85}), while the slow pulse is unstable (Flores~\cite{Flores}, Evans~\cite{EvansIV}, Ikeda, Mimura and Tsujikawa~\cite{IMT87}). In addition, there are stable fast pulses with oscillatory tails (Carter and Sandstede~\cite{CS15}, Carter, de Rijk and Sandstede~\cite{CdRS16}).

By above, each pulse $\Phi$ on $\cS=\R$ defines a smooth {\it axisymmetric} traveling wave solution $\bu(x, \theta, t)=\Phi(x-ct)$ of Eq.~(\ref{FHN-S1}$\cS_R$), propagating along the cylindrical axis. Its speed $c$ is determined by
the parameters $\alpha$, $\gamma$, and~$\eps$. After \cite{TBS}, we continue to call them the pulses. We are interested in the fast pulse, for which we keep the notation $\Phi(x-ct)$. 

\paragraph{\bfseries Results of \cite{TBS}.} As was mentioned above, an analysis of the FHN system (\ref{FHN-S1}$\cS$) on cylindrical surfaces was initiated in \cite{TBS}. It is shown in \cite{TBS} that, for $\cS = \cS_R, $ with $R\leq 1, $ mild solutions which are initially close to $\Phi(x-ct)$ approach nearby translates of~$\Phi(x-ct)$ as $t\to \infty$, or put differently,

\begin{enumerate}[(i)]

\item On a cylinder of small constant radius, the (fast) pulses are  asymptotically orbitally stable under general perturbations of the initial conditions that depend on both spatial variables.
\end{enumerate}

One does not expect that there are pulse-like traveling waves on generic cylindrical surfaces, $\cS_\rho$, even if $\rho$  is close to a constant $R\leq 1$. However, it is shown in \cite{TBS} that:

\begin{enumerate}[(ii)]

\item Eq.~(\ref{FHN-S1}$\cS_\rho$), for the warped cylinder $\cS_\rho, $ has stable, nearly radial pulses surrounded by small fluctuations propagating 
along the cylindrical axis (as is the case with real axons).
\end{enumerate}

Put differently, on a warped cylinder, solutions that are initially close to a pulse stay near a propagating pulse for all time. (Note that the restriction $ R \leq 1 $ is not optimal. In particular, we conjecture that the results of \cite{TBS} hold for $R<1/\sqrt{1-\al}. $ However, the optimal interval of $R$ for which the results of \cite{TBS} hold is not known. We thank an anonymous referee for posing this question.)

In this paper, we consider the {\it undulated} cylinders, i.e.\ graphs, \eqref{graph-cyl}, over standard cylinders with the defining radius function, $\rho, $ radially symmetric, i.e.\ depending only on $x$, $\rho=\rho(x)$. We prove:
\begin{itemize}
\item A spontaneous symmetrization of solutions, initially close to radially symmetric ones;
\item Approximation of solutions to (\ref{FHN-S1}$\cS_\rho$) by solutions of a system in 1 spatial dimension (Eq.~\eqref{eq:radlapl} below) with coefficients depending on the radius profile $\rho=\rho(x).$
\end{itemize}

To formulate our results precisely, we let $\vec L^2 \equiv  L^2(\cS_\rho, \R^2)$ be the Hilbert space of square integrable functions on $\cS_\rho $ (see \eqref{vecL2-rho} below for more a detailed description) and define the vector Sobolev space, $ \vec{H}^{1,0}  \equiv   H^{1,0}(\cS_\rho, \R^2), $ in a standard way as
\begin{equation}\label{eq:Hkl}
\vec{H}^{1,0}
:=\left\{\bu\in \vec L^2
\ \big\vert\  
(-\Delta)^{1/2} u_1 \in L^2\right\}\,,
\end{equation}
where $ \Delta\equiv \Delta_{\cS_\rho}$ is the Laplace-Beltrami operator on $ \cS_\rho $ (see \eqref{laplaxial} below) and $L^2 \equiv L^2(\cS_\rho, \R)$, with norm
\begin{equation}\label{normk0}
\|\bu\|_{1,0}^2:= \|(-\Delta)^{1/2}u_1\|^2 + 
\eps^{-1}\|u_2\|^2\,,
\end{equation}
and its dual $ \vec{H}^{-1,0}.$ Furthermore, we say (\ref{FHN-S1}$\cS_\rho$) has a {\it strong}, local solution $\,\bu(t) \, $ in $ \, \vec{H}^{1,0} , \, $ if 
$$\bu \in C([0,T],\vec{H}^{1,0}) \cap C^1([0,T], \vec{H}^{-1,0}) \,,$$ 
for some $ T>0, $ and $\,\bu(t)\,$ satisfies Eq.~(\ref{FHN-S1}$\cS_\rho$). We say that the solution $ \bu $ is {\it global} if we can take $ T=\infty.$ 

Consider (\ref{FHN-S1}$\cS_\rho$) with $\rho$ of class $C^2$, bounded and radially symmetric (i.e.\ independent of $\theta$) with $ \|\partial_x\rho\|_{L^\infty}<\infty. $ We consider strong solutions of this equation on an interval $ [0,T], \; T>0, $ satisfying
\begin{equation}\label{e25225}
\sup_{0\leq t\leq T} \|\bar\bu(t)\|_{1,0} \leq K , 
\end{equation}
where $\bar f(x):=\frac1{2\pi}\int_0^{2\pi}f(x, \theta) d \theta$, and denote the class of such solutions by $ S_{T,K}.$

A local existence of (\ref{FHN-S1}$\cS_\rho$) on $[0,T] $ with $T$ independent of the parameters $\al, \gamma $ and $\eps $ is standard and is discussed in Appendix \ref{sec:locglobexist}. It is proven in exactly the same way as a higher regularity local existence is proven in \cite{TBS}, Proposition 2.1. 
One can readily show that the local solutions satisfy \eqref{e25225} with constants $K $ independent of $\al, \gamma $ and $\eps, $ so that the sets $ S_{T,K} $ are definitely non-empty.

Below, we prove a priori bounds on solutions of (\ref{FHN-S1}$\cS_\rho$) independent of the existence results.

\allowdisplaybreaks
Our first result is given in:

\begin{theorem}
\label{thm:rad-conv}
There exist $ r_* $ and $ \del_0 $ s.t., for any $0 < \rho \le r_* $ and any $T>0$ and $K\geq 1, $ any strong solution $\bu(t)$ to (\ref{FHN-S1}$\cS_\rho$) in $S_{T,K}, $\, with an initial condition $  \bu_0 \in \vec{H}^{1,0} $ obeying $ \, \|\bu_0-\bar{\bu}_0\|_{1,0}\leq \del_0/K $, satisfies the estimate
\begin{align}\label{rad-conv}
\|\bu(t)- \bar\bu(t) \|_{1,0} 
\; \le \; C \;e^{-\tfrac{\gamma}{2}\,\eps \,t}\; \|\bu_0- \bar\bu_0 \|_{1,0} \qquad (0\leq t\leq T),
\end{align}
where $C$ is a positive constant independent of $ T, K $ and $ \eps$.
\end{theorem}

A proof of Theorem \ref{thm:rad-conv} is given in Section \ref{sec:rad-collapse} with technical estimates carried out in Section \ref{sec:apriori}.

\begin{remarks}\label{r27225} (i)  We expect that for initial conditions in Theorem \ref{thm:rad-conv}, the solutions exist and belong to $S_{T,K}, $  $ K=\mathcal O(1) $ and $ T=\infty. $ Let $$\vec H^{2,1}
:=\left\{\bu\in \vec L^2 :  
\Delta u_1 \in L^2\,, \; \p_x u_2\in L^2 \right\}\;,$$ 
(cf.~\eqref{gvSs} in Appendix \ref{sec:locglobexist}) with the norm  $ \|\bu\|_{2,1}= (\|(\Delta u_1, \partial_xu_2)\|^2 + \|\bu\|^2)^{1/2}\,,$ and define the manifold of pulses
$$\mathcal M := \{\Phi_h|\, h\in\R\},$$
where, recall, $\Phi_h(x)=\Phi(x-h). $ Let also $\operatorname{dist}(\bu, \mathcal M)=\inf_{\bm{v}\in\mathcal M}\|\bu-\bm{v}\|_{2,1}. $
It is shown in \cite[Theorem 1.2]{TBS} that there are $\kappa, \delta_*>0 $ such that if $ R^{-1}\|\rho-R\|_{C^2}\leq \delta_* $ for some $R\leq 1, $ then for any initial condition $  \bu_0 $ s.t.\ $\operatorname{dist}(\bu_0, \mathcal M)\leq \kappa, $ there exists the unique global solution $\bu(t)$ of (\ref{FHN-S1}$\cS_\rho$) satisfying $ \operatorname{dist}(\bu(t), \mathcal M)\leq C (\kappa+\delta_*), $ for some constant $ M, $ independent of $\epsilon.$ The latter bound implies that $ \|\bu\|_{1,0} \leq  \|\Phi\|_{H^1} + C (\kappa+\delta_*)  = \mathcal O(1) $ and therefore $\bu(t)\in S_{T,K} $ with $T=\infty $ and $R=\mathcal O(1).$

(ii) Unlike in the results of \cite{TBS}, the surfaces $\cS_\rho $ in Theorem \ref{thm:rad-conv} are not assumed to be closed to a straight cylinder.

(iii) The results above depend strongly on the sectional curvatures $1/\rho(x)$ of the sections of $\cS_\rho $ by planes orthogonal to the $x$-axis. In particular, $r_*$ in Theorem \ref{thm:rad-conv} is the lower bound on $\rho(x).$ It is crucial here that the new $\theta$-degree of freedom is compact.

Surprisingly, the result above is independent of the principal curvature along $x$-axis and depends only on the line element $\sqrt{1+(\partial_x\rho)^2}dx$ along $x. $ 

We conjecture that our result extend to warped surfaces, i.e.\ to $\rho$ dependent not only on $x$ but also on $\theta.$
\end{remarks}

We also introduce the radial  $L^2$ and Sobolev spaces in 1 dimension
\begin{flalign}\label{L2rhorad}
&\vec L^2_{rad} \equiv L^2(\mathbb R,\sqrt{g}dx, \mathbb R^2)\,, \quad  L^2_{rad} \equiv L^2(\mathbb R,\sqrt{g}dx, \mathbb R)\,, \\[0.5em]
&\vec{H}_{rad}^{1,0} \;:=\; \left\{\bu\in \vec L^2_{rad}: (-\Dr)^{1/2}u_1\in L^2_{rad}\right\}\,,
\label{Hk0rhorad}
\end{flalign}
where $\Dr $ is the radial part of the Laplace-Beltrami operator $\Delta, $ see \eqref{Laplrad} and \eqref{laplaxial} below, with the norms defined as in \eqref{normk0} (actually, they could be defined by \eqref{normk0}) and denoted by the same symbols and similarly for the scalar products. Moreover, we also use the dual space $ \vec{H}_{rad}^{-1,0} $ to $\vec{H}_{rad}^{1,0}.$

The next result compares the averages $ \bar\bu(t) $ of solutions $ \bu(t) $ to (\ref{FHN-S1}$\cS_\rho$) to solutions $\,\bm{w}=(w_1,w_2)\, $ of the following system of equations
\begin{equation}\label{eq:radlapl}
\begin{split}
\partial_t w_1 & \;= \; \Dr \,w_1 \,+\, f(w_1)\, -\, w_2 \,,
\\[0.5em]
\partial_t w_2 & \; = \; \eps \, (w_1 \,-\,\gamma\, w_2)\,.
\end{split}
\end{equation}
We consider strong solutions of Eqs~(\ref{FHN-S1}$\cS_\rho$) and \eqref{eq:radlapl} in the spaces
\begin{eqnarray*}
V_1 &:=& C([0,T],\vec{H}^{1,0}) \cap C^1([0,T],\vec{H}^{-1,0}) \;,
\\
V_{1,rad} &:=& C([0,T],H^{1,0}_{rad}) \cap C^1([0,T],H^{-1,0}_{rad})\;.
\end{eqnarray*}
We have:
\begin{theorem}\label{thm:average} Under the conditions of Theorem \ref{thm:rad-conv}, the average $\,\bar\bu(t)\, $ of the strong solution $ \, \bu(t)\,$ of (\ref{FHN-S1}$\cS_\rho$) in $ V_1 $ and the strong solution $\,\bm{w}\, $ in $\,  V_{1,rad}\,$ to \eqref{eq:radlapl}, with initial conditions $\bu_0$ and $\bm{w}_0$, respectively, satisfy the estimate

\begin{equation}\label{vest}
\|\bar\bu(t)- \bm{w}(t) \|_{1,0}
\;\leq \;C \,  e^{\mu \, t}\,\left(\|\bar\bu_0-\bm{w}_0\|_{1,0}+\|\bu_0-\bar\bu_0\|_{1,0}^2\right)\;,
\end{equation}
for some constants $ C,\mu>0 $ independent of $ T, R $ and $\eps $ and $ \, 
t$ $\ll$ $-\ln\big(\|\bar\bu_0-\bm{w}_0\|_{1,0}$ $+\|\bu_0-\bar\bu_0\|_{1,0}^2\big).$  
\end{theorem}

\noindent A proof of Theorem \ref{thm:average} is given in Section \ref{seqefeq}. Eq.~\eqref{eq:radlapl} is a system in 1 spatial dimension depending on the radius profile $\rho=\rho(x). $ It is considerably simpler for mathematical analysis and numerical simulations, and will be explored in a future study. It offers an \textit{effective system} for the pulse propagation.

\begin{remarks}\label{2r27225}
\vspace{0.6em}

(i) To keep exposition simple, we restricted ourselves to the energy space $\vec{H}^{1,0}. $ It would be interesting to extend the results above to higher regularity. It is natural to use Sobolev spaces with different smoothness in different components, since (\ref{FHN-S1}$\cS_\rho$) is of the second order in $u_1 $ and the zero order in $u_2.$

(ii) The global existence in the space $ \, \vec{H}^{2,1},\, $ defined in Remark \ref{r27225}(i) and Appendix \ref{sec:locglobexist}, for certain class of initial data is proven in \cite{TBS}. One might be able to extend the proof of \cite{TBS} to $ \, \vec{H}^{1,0} \,$ under conditions similar to those of \cite{TBS}.

(iii) There are a few results in higher dimensions. In $\R^2$, 
Mikhailov and Krinskii~\cite{MK83}, and Keener~\cite{Keener86} 
studied spiral solutions of Eq.~(\ref{FHN-S1}$\R$). In $N$-dimensions,  for $N \geq 2$, Tsujikawa, Nagai, Mimura, Kobayashi and Ikeda~\cite{Tsujikawa_89} proved that there exist stable fast pulse solutions propagating in a one-dimensional direction. Numerical computation of solitary waves in cylindrical domains was done in \cite{LPSS} (see \cite{Sandst} for a review and further references).

(iv) Existence and stability for fast pulses have been studied for variants of Eq.~(\ref{FHN-S1}$\R$), where the second equation also has a diffusion term (Cornwell and Jones~\cite{Cornwell-Jones18}, Chen and Choi~\cite{Chen-Choi15}, Chen and Hu~\cite{Chen-Hu14}). It would be interesting to extend the results above to undulated surfaces variants of these equations.

(v) Another system that admits stable fast pulses is the 
discrete analogue of Eq.~(\ref{FHN-S1}$\R$) (Hupkes and Sandstede~\cite{HS13}, 
Schouten-Straatman and Hupkes~\cite{SSH19}, Hupkes, 
Morelli, Schouten-Straatman and Van Vleck~\cite{HMSSVV}).

(vi) Another natural extension of Eq.~\eqref{FHN-S} would be to have an $x$-dependent nonlinearity modeling varying mylination.
\end{remarks}

\vspace{0.6em}

\paragraph{\bfseries Outline of the approach.} Our main tool is the method of differential inequalities for quadratic functionals (cf.\ \cite{AKS,FGJS,KP13} and references therein). Namely, combining the functionals
\begin{equation}
X_k(\bu^\perp) :=\frac12 \langle  \bu^\perp, (-A)^k\bu^\perp \rangle \; ,
\end{equation}
with $k=0,1$, where $\bu^\perp=\bu-\bar\bu, $  
and
\begin{equation}\label{e21225}
A:=\begin{pmatrix} \Del -\alpha & -1 \\
\eps & - \eps\gamma 
\end{pmatrix}\,,
\end{equation}
where $\Del \equiv \Delta_{\cS_{\rho}},$  into a simple functional $X(\bu^\perp), $ we obtain the differential inequality, for $X(t)\equiv X\big(\bu^\perp(t)\big), $ where $\bu(t) $ is a solution of (\ref{FHN-S1}$\cS_\rho$) in $ S_{T,K} \,,$
\begin{equation*}\label{blan3}
\partial_t X \leq - \tfrac{\gamma}{2}\,\eps \, X \,+\,  C_3 \,K^2 \,X^2 
\,+\,C_4 \,X^3\,.
\end{equation*}
Integrating this inequality (see abstract Lemma \ref{lemdifineq}) and using that $ X(t) \simeq \|\bu^\perp(t)\|^2_{1,0}, $ we arrive at Eq.\ \eqref{rad-conv}.

Theorem \ref{thm:average} is proven similarly, except that the coefficient in front of an equivalent of $X(t) $ is positive which gives a time growing exponential. The later result could possibly be improved but this would require additional powerful machinery, similar to that used in \cite{TBS}, such as detailed spectral analysis of the linearized system and comparison estimates for solutions of (\ref{FHN-S1}$\cS_\rho$) and (\ref{FHN-S1}$\cS_R$).

The paper is organized as follows. In Section \ref{sec:rad-collapse}, we give a proof of Theorem \ref{thm:rad-conv} with technical estimates carried in Section \ref{sec:apriori}. In Section \ref{seqefeq}, we prove Theorem \ref{thm:average} and in Appendix \ref{sec:locglobexist}, we sketch a proof of the local existence for system (\ref{FHN-S1}$\cS_\rho$).

\noindent
We conclude the introduction with a discussion of the spaces used in this work.   

\vspace{0.6em}

\paragraph{\bfseries Remarks on spaces and operators.} 
In the cylindrical coordinates, undulated cylindrical surface $\mathcal S_\rho,\; \rho=\rho(x), $ cf.~\eqref{graph-cyl}, is written as
\begin{align} 
\label{Srho}
\cS_\rho :
&= \left \{(x,\rho(x)\cos\theta, \rho(x)\sin\theta)
\in\R^3 \ 
\big\vert\ 
x\in \R, \theta\in [0,2\pi) \right\},
\end{align}
and the area element on $\cS_\rho$ is given by $d\mu=\sqrt{g} dx d \theta$, where $g$ is the determinant of the metric tensor, defined in \eqref{Laplrad},
\begin{align}\label{g-rad}
&g=(1+\rho_x^2) \rho^2 \,,\qquad \rho_x=\p_x\rho.
\end{align}

In the cylindrical coordinates, the space $\vec L^2 \equiv L^2(\cS_\rho, \R^2) $ is identified with   
\begin{equation}\label{vecL2-rho}\vec L^2 \equiv L^2(\cS_\rho, \R^2)  \equiv L^2(U, d\mu, \R^2) ,\end{equation}
where, recall, $U:=\R\times (\mathbb R\setminus (2\pi\mathbb Z))$ and $d\mu:=\sqrt{g} dx d \theta$. We define the  inner product on $\vec L^2$ as 
\begin{equation}\label{inner-rho}
\langle \bm{u},\bm{w} \rangle :=\int_{U} 
(u_1 w_1+\eps^{-1}u_2 w_2)\, d\mu\,.
\end{equation}
The corresponding norm will be denoted by $\| \cdot \|$. The scaling in \eqref{inner-rho} was introduced in \cite{TBS} in order to make the linearized operator \eqref{e21225} symmetric.

Turning to the Laplace-Beltrami operator $ \Delta\equiv \Delta_{\cS_\rho}$ on $ \cS_\rho, $ with $ \,\rho=\rho(x) \, $ independent of $\,\theta, $ which we use to define the Sobolev spaces, in the cylindrical coordinates, it is given by
\begin{equation}\label{laplaxial}
\Del
\,\equiv\,
\Delta_{\cS_{\rho}}
\,=\,
\frac{1}{\sqrt{g}} \big[ \p_x \frac{\rho^2}{\sqrt{g}} \,\p_x  \,+\,  \p_{\theta} \frac{1+\rho_x^2}{\sqrt{g}} \,\p_{\theta} \big]\;,
\end{equation}
where, recall, $\rho_x=\p_x\rho, $ and $g$ is given in \eqref{g-rad}. The operator $\Delta$ is self-adjoint and non-positive on the space $ L^2 \equiv L^2(\cS_\rho, \R) $ of scalar $L^2$-functions on $\cS_\rho $ (cf. \eqref{vecL2-rho}).

The radial part $\Dr $ of the Laplace-Beltrami operator $\Delta \equiv \Delta_{\cS_\rho} $ is given by
\begin{equation}\label{Laplrad}
\Dr \; = \; \frac{1}{\sqrt{g}} \; \p_x \frac{\rho^2}{\sqrt{g}} \,\p_x \;,
\qquad
g=\rho^2(1+ \rho_x^2)\,.
\end{equation}

\paragraph{\bfseries Notation.} To summarize, we denote the $L^2$ and Sobolev spaces as $L^2 \equiv L^2(\cS_\rho, \R) $ and 
$$H^1 \equiv H^1(\cS_\rho,\mathbb R) :=\left\{f\in L^2
\ \big\vert\  
(-\Delta)^{1/2} f \in L^2\right\}\,,$$
for scalar functions, and as $ \vec L^2 \equiv  L^2(\cS_\rho, \R^2)$ and $\vec H^{1,0} \equiv H^{1,0}(\cS_\rho, \R^2), $ for the vector functions.

In what follows, the norms in $L^2 \equiv L^2(\cS_\rho, \R), \, \vec L^2 \equiv  L^2(\cS_\rho, \R^2) , \, H^1 \equiv H^1(\cS_\rho, \R)$ and $\vec H^{1,0} \equiv  H^{1,0}(\cS_\rho, \R^2)$, are denoted by $\|\cdot\|,\, \|\cdot\| ,\, \|\cdot\|_{H^1} $ and $ \|\cdot\|_{1,0}, $ respectively. The function $ \rho(x) $ is fixed and we omit it in the notation, except for the notation $\cS_\rho$ for the cylindrical surface.

Furthermore, the parameters $ \al, \gamma $ and $ \eps $ and the function $\rho $ entering (\ref{FHN-S1}$\cS_\rho$), \eqref{nonlinf} are fixed positive numbers. $ \eps $ is essential in the stability theory and we trace the dependence on it while ignoring $ \al $ and $ \gamma. \; C, C_1, $ etc are arbitrary constants independent of $ T, K $ and $ \rho. $ These constants may depend on $\al$ and $\gamma$ and the Sobolev constants and may vary from one inequality to another. The maximum of various Sobolev constants in 1D and 2D is denoted by $M. $ The constants used in Sections \ref{sec:rad-collapse} and \ref{seqefeq} are not related to each other.

Finally, for two positive functions $f$ and $g, $ the notation $f\simeq g$ say that $ c_1 f \leq g \leq c_2 f $ for some positive constants, which might depend on $\al $ and $\gamma $ but are independent on $R, T , \eps $ and $\rho.$

We use the following typographical convention. Constants are denoted by lower greek and capital roman letters (e.g.\ $\al, \gamma, \eps, T, K, C, $ etc), operators, by roman capitals (e.g.~$A$), functions by either lower or capital roman letters (e.g.~$f, g, h $ or $W, X, X_k, Y, Z, $ to eliminate the confusion in the latter case, we display the arguments, i.e.~$W(t), X(t),\, X_k(t), \, Y(t) $ and $ Z(t)$), and maps, by capital math calligraphic letters (e.g.\ $\mcF$ and $\mcN$). The only exception is the function $ \rho(x).$ Vector functions are denoted by boldface letters, e.g. $\bm{f}=(f_1,f_2).$

\section{Angular contraction: Proof of Theorem~\ref{thm:rad-conv}}\label{sec:rad-collapse}

\subsection{Proof of Theorem \ref{thm:rad-conv} modulo estimates (\ref{dotx0est}) and~(\ref{dotx1est})}\label{sec28225}

In this subsection we prove Theorem \ref{thm:rad-conv} modulo two key technical estimates demonstrated in the next subsection. We split the proof of Theorem~\ref{thm:rad-conv} into several steps.

To begin with, we observe that on the cylindrical surface $\cS_\rho$, the FHNcyl system, (\ref{FHN-S1}$\cS_\rho$), takes the form
\begin{equation}\label{FHNrho1}
\begin{split}
\p_t u_1 &= \Delta u_1 + f(u_1) - u_2 \,,\\
\p_t u_2 &= \eps (u_1 - \gamma u_2)\,,
\end{split}
\end{equation}
where $\Delta \equiv \Delta_{\cS_{\rho}}, $ the Laplace-Beltrami operator on $ \mathcal S_\rho $ (see \eqref{laplaxial} for an explicit expression).

\vspace{0.6em}

\paragraph{\bfseries Canonical form of equation \textbf{(\ref{FHNrho1})}.} Denote the right hand side of Eq.~\eqref{FHNrho1}
by $\mcF(\bu)$. 
Then the  initial-value problem for \eqref{FHNrho1} becomes
\begin{equation}\label{ivp-Frho}
\p_t \bu = \mcF(\bu), \qquad \bu\vert_{t=0} = \bm{u_0}.
\end{equation}
We consider this problem on the space $\vec{H}^{2,0}$. We rewrite the initial-value problem Eq.~\eqref{ivp-Frho} by expanding $\mcF(\bu)$ in $\bu$ and using $\mcF(0)=0$ to obtain 
\begin{equation}\label{Frho}
\mcF(\bu) = A \bu + \mcN(\bu) \;,
\end{equation}
where $ \, A \bu =  d\mcF(0)\bu \, $ and $ \, \mcN(\bu)  \, $ is defined by this equation, namely, $\mcN(\bu):=\mcF(\bu)-A \bu$. Now, Eq.~\eqref{ivp-Frho} becomes
\begin{equation}\label{rho-eq-canon}\p_t \bu = A \bu + \mcN(\bu) \; ,\end{equation}
where $ A $ is the principal part of $\mcF$ given by its G\^ateaux derivative at $ \bu=0$

\begin{equation}\label{A} 
A:= d\mcF(0)
= \begin{pmatrix} \Delta -\alpha & -1 \\
\eps & - \eps\gamma 
\end{pmatrix}\,,
\end{equation}
with $ \Delta $ the Laplace-Beltrami operator given in \eqref{laplaxial}, and $\mcN(\bu)$ is the nonlinear part (nonlinearity) given by
\begin{equation}
\label{N}
\mcN(\bu) := \mcF (\bu) - A \bu = 
\begin{pmatrix}
h(u_1) \\
0
\end{pmatrix}, \quad h(u_1):=-u_1^3 +(\alpha+1)u_1^2 \;.
\end{equation}

\vspace{0.6em}

\paragraph{\bfseries Operator $\bm{A}$.} Recall the space $ \vec{H}^{1,0} \equiv  H^{1,0}(\cS_\rho, \R^2) $ and define the space  $ \vec{H}^{k,0}$ as
\begin{equation}\label{e27225}
\vec{H}^{k,0}
:=\left\{\bu\in \vec L^2
\ \big\vert\  
(-\Delta)^{k/2} u_1 \in L^2 \;\;\mbox{and}\; u_2\in  L^2 \right\}\,,
\end{equation}
$k=1,2. $ By the standard elliptic theory, the G\^ateaux derivative  $A=d\mcF(0)$ maps $ \vec H^{k,0} $ into $  \vec H^{k-2,0}, $ for $k=1,2 $ (see \cite{TBS}). 

\vspace{0.6em}

\paragraph{\bfseries Lyapunov-Schmidt decomposition.}
Let  $P$ be the projection on the average of the function in $\theta$ (zero harmonic):
\begin{equation}\label{proj}
(Pf)(x) := \bar f(x)  \equiv \frac1{2\pi}\int_0^{2\pi}f(x, \theta) d \theta ,
\end{equation}
and similarly for vector-functions, and let $P^\perp:=\one-P$. Note that the operators $ P, P^\perp $ commute with $ \Delta $ (on scalar functions) and with $A$ (on vector-functions). Furthermore, let
\begin{equation}\label{decomp}
u^\perp:=P^\perp u=u-\bar u \, , 
\end{equation}
\begin{equation}\label{Aperp}
A^\perp:= P^\perp A 
\; = A P^\perp 
\;
= P^\perp A P^\perp
\qquad\mbox{and}\qquad
\mcN^\perp(\bu) := P^\perp \mcN(\bu) \,.
\end{equation}

\noindent Note that, since $\rho$ is radially symmetric, the operator $A$ commutes with the rotations and therefore with $P$ and $P^\perp$.

By the definitions above, if $ \bu(t) $ is a solution to \eqref{rho-eq-canon}, then $ \bu^\perp(t) := P^\perp\bu(t) $ satisfies the equation
\begin{equation}\label{duperp}
\partial_t\bu^\perp = A \bu^\perp + \mcN^\perp(\bu) \,,
\end{equation}
obtained by applying $ P^\perp $ to Eq.~\eqref{rho-eq-canon} and using that 
$ P^\perp $ and $ A $ and $ P^\perp $ and $ \partial_t $ commute. (Note that, since $u^\perp=P^\perp u^\perp, $ we have $A \bu^\perp=A^\perp \bu$.)

\vspace{0.4em}

\paragraph{\bfseries Lyapunov functions} To estimate various norms of $ \bu^\perp(t), $ we use differential inequalities for the following Lyapunov functions
\begin{equation}\label{defX12}
X_k(t)\equiv X_k\big(\bu^\perp(t)\big):=\frac12 \langle  \bu^\perp(t), (-A)^k\bu^\perp(t) \rangle \; ,
\end{equation}
with $k=0,1$. Note that $X_0(t)=\tfrac12\|\bu^\perp(t)\|^2$ and, as shown below, that $X_1(t) \simeq \|\bu^\perp(t)\|^2
_{1,0}.$

Recall that $\Del $ is defined in \eqref{laplaxial}, and let $\nabla\equiv \nabla_\rho $ be the corresponding gradient. Using the formula \eqref{laplaxial} and then integrating by parts, we find
\begin{equation}\label{2e26225}
\langle f^\perp , \Delta f^\perp\rangle = -
\left\| \nabla f^\perp \right\|^2 =  -\int\left(\tfrac{\rho^2}{g}|\partial_x f|^2 + \tfrac{1+\rho_x^2}{g}|\partial_\theta f|^2\right)\sqrt{g}\,dxd\theta.
\end{equation}
Using the definition of $A$ in \eqref{e21225} and the definition of the inner product in \eqref{inner-rho} and using the relation \eqref{2e26225}, we obtain
\begin{eqnarray}\label{def-x1}
2 X_1(t) 
& = &
- \langle  \bu^\perp(t), A\bu^\perp(t) \rangle
\\[0.5em]
&=&
-\langle u_1^\perp , \Delta u_1^\perp\rangle
+\al \|u_1^\perp\|^2
+\gamma \|u_2^\perp\|^2
\nonumber
\\[0.5em]
&=&
\| \nabla u_1^\perp \|^2
+\al \|u_1^\perp\|^2
+\gamma \|u_2^\perp\|^2
\label{def-x1b}
\\[0.5em]
&\simeq &
\|\bu^\perp(t)\|^2
_{1,0} \,,
\label{def-x1c}
\end{eqnarray}
where recall the notation $\, g \simeq f\, $  says that $ f $ and $ g $ satisfy $  \,c_1 f \leq  g \leq c_2 f\, $ for some positive constants $\, c_1, c_2 \,$ depending on the fixed parameters $ \al,\gamma $ and $ \eps $ of Eq.~(\ref{FHN-S1}$\cS_\rho$).

The point here is that on $\, \operatorname{Ran}(P^\perp) \, $ the operator $\, -A \, $ is coercive: for any $r_*>0,$
\begin{align}
\label{Arho-perp-est} \langle A^\perp \bu^\perp, \bu^\perp \rangle
&\le -\,\frac12 \, \|\n u_1^\perp\|^2 \,-\,\frac{c}{2r_*^2}\, \|u_1^\perp\|^2
\,-\,\gamma\, \|u_2^\perp\|^2 \;,
\end{align}
with $ c = (1+\|\rho_x\|^2_{L^\infty}) $, provided $ \rho\leq r_*. $ Indeed, using \eqref{A} and the inner product \eqref{inner-rho}, we compute
\begin{align}
\label{A-perp-ip}   \langle A^\perp \bu^\perp, \bu^\perp \rangle
&= \langle (\Delta^\perp -\alpha) u_1^\perp, u_1^\perp\rangle -\gamma \|u_2^\perp\|^2,\end{align}
where $\Delta^\perp $ is defined according to \eqref{Aperp}. Writing $\langle \Delta^\perp  u_1^\perp, u_1^\perp\rangle=\frac12\langle \Delta^\perp  u_1^\perp, u_1^\perp\rangle  + \frac12\langle \Delta^\perp  u_1^\perp, u_1^\perp\rangle $ and using that $\al>0 $ and the inequality \eqref{Lapl-est} given below, we find \eqref{Arho-perp-est}.

\begin{lemma}\label{lem:Lapl-est}
Let $ \rho\leq r_*. $ Then we have the estimate
\begin{align}\label{Lapl-est} \langle \Delta^\perp  u_1^\perp, u_1^\perp\rangle \le - \frac{c}{r_*^2}\, \|u_1^\perp\|^2\,,\qquad c = 1+\|\rho_x\|^2_{L^\infty}\,.\end{align}
\end{lemma}
\begin{proof} Using \eqref{2e26225} and the definition \eqref{g-rad} giving $ g=(1+\rho_x^2) \rho^2\leq (1+\rho_x^2) r_*^2, $ provided $\rho\leq r_*, $ we find
\begin{equation}\label{e26225}
\langle u_1^\perp , \Delta u_1^\perp\rangle
\leq
-\frac{c}{r_*^2}
\int |\partial_\theta u_1^\perp|^2\sqrt{g}\,dx d\theta 
=
\frac{c}{r_*^2} \langle u_1^\perp , \partial_\theta^2 u_1^\perp\rangle \,.
\end{equation}
On the subspace $ \operatorname{Ran}(P^\perp), $ the spectrum of the operator $\partial_\theta^2 $ is $\{-n^2, n=1,2,\ldots\}. $ Hence, by a standard spectral theory, $\langle u_1^\perp , \partial_\theta^2 u_1^\perp\rangle \le -\|u_1^\perp\|^2, $ which together with \eqref{e26225} implies \eqref{Lapl-est}.
\end{proof} 

In the next subsection, assuming $\rho\leq r_* $, we prove the following estimates (see Propositions \ref{prop:domain-A} and \ref{propdX1}): for any $ T>0,\, K \geq 1,$ and $ \bu\in S_{T,K}\,, $
\begin{equation}
\partial_t X_0 \le - \frac{1}{4} \, \|\n u_1^\perp\|^2 \;-\; \frac{c}{4 \, r_*^2} \, \|u_1^\perp\|^2 \;-\; \gamma\, \|u_2^\perp\|^2
\; + \; C_1 \,K^2 \, \|u_1^\perp\|^4 \; ,
\label{dotx0est}
\end{equation}
and
\begin{equation}
\label{dotx1est}
\begin{split}
\partial_t X_1 
\, \leq& 
\; - 
\, \frac{1}{8}\,
\|\Delta u_1^\perp\|^2
\;-\;
\bigg(\frac{1}{16} \bigg(\frac{c}{r_*^2}\bigg)^2- C_3 K^4\bigg)\; \|u_1^\perp\|^2 
\\
&\, + \,C_2 \big( \|u_2^\perp\|^2
\,+\,K^2\,X_1^2 
\,+ \, X_1^3 \big)
\;,
\end{split}
\end{equation}

\noindent for some positive constants $\, C_2 $ and $ C_3 $ independent of the parameters $ T, K $ and $ \eps$ and of the function $\rho $ and satisfying 
\begin{equation}\label{e28225}
C_2=\frac18 \big(11 - 5\eps\gamma^2\big).
\end{equation}

Note that \eqref{dotx0est} implies a bound on $ X_0(t)=\tfrac12\|\bu^\perp(t)\|^2. $ For this, we use the following general result:

\begin{lemma}\label{lemdifineq}
Let $ \, X(t)\, $ be a positive, differentiable function satisfying differential inequality 
\begin{equation}\label{gdifeq}
\partial_t X \leq - \nu \, X \,+\, M_1 \,X^2 \,+\,M_2 \,X^3 \,,
\end{equation}
with an initial condition $ \, X(0)\ll \min\left(1/M_1, 1/\sqrt{M_2}\right). $ Then $ \, X(t) \, $ satisfies the estimate 
\begin{equation}\label{lemdifineq1}
X(t) \, \leq\,   2 \, X(0) \,   e^{-\nu\,t}\,, \qquad t\geq 0.
\end{equation}
\end{lemma}
\begin{proof}Let \, $ Y(t) := e^{\nu t} X(t). $ Then using \eqref{gdifeq}, we obtain
\begin{equation}\label{dtyest1}
\partial_t Y(t) \, \leq \, M_1 \,e^{-\nu t} \, Y^2(t) \, + \, M_2 \, e^{-2\nu t} \, Y^3(t)\,.
\end{equation}
Integrating this inequality, and taking the supremum of the result over the interval $ [0,s] $ and using the notation
\begin{equation}\label{dtyest2}
Z(s) :=\sup\limits_{0\leq t\leq s} Y(t)\,,
\end{equation} 
we obtain the inequality
\begin{equation}\label{dtyest3}
Z(s) \leq X(0) \,+\, M_1\, Z^2(s) \,+ \,M_2\,Z^3(s)\;.
\end{equation}
For $ X(0) \ll  \min\left(1/M_1, 1/\sqrt{M_2}\right) , $ the polynomial  $\, p(Z) = M_2 Z^3 +M_1 Z^2 - Z + X(0)\, $ has three roots (see Fig.~\ref{fig:pol})
\begin{equation*}\label{roots} X(0)<Z_1 <2 X(0)\,,\qquad Z_2\approx 1 \,, \qquad Z_3 < 0 \, .\end{equation*}

\begin{figure}[H]
\begin{center}

\begin{tikzpicture}[domain=-2:4,xscale=2.7,scale=0.7]
\draw[<->,thick] (0,3.5) node[right]{$p(Z)$} -- (0,0) -- (4,0) node[below] {\footnotesize$Z$};

\draw[thick] (0,0)--(-3,0) ;
\draw[thick] (0,0)--(0,-2) ;


\draw[mark size=1.3pt] (-2,0) node[] {\pgfuseplotmark{*}};
\draw[] (-2,0) node[below right=-0.2em]  {\scriptsize $Z_3$};

\draw[mark size=1.3pt, black] (3,0) node[] {\pgfuseplotmark{*}};
\draw[] (3,0) node[below right=-0.2em]  {\scriptsize $Z_2$};

\draw[mark size=1.3pt, black] (1,0) node[] {\pgfuseplotmark{*}};
\draw[] (1,0) node[below left=-0.2em]  {\scriptsize $Z_1$};

\draw[thick, dashed, domain=-2.3:3.3] plot ( {\x},  { 0.2*(\x - 3)*(\x - 1)*(\x + 2) } );


\draw[mark size=1.3pt, black] (0,1.2) node[] {\pgfuseplotmark{*}};
\draw[] (0,1.2) node[above right=-0.2em]  {\scriptsize $X(0)$};

\draw[black] (0,0) node[below left=-0.2em] {\scriptsize$0$};
\draw[mark size=1.5pt, black] (0,0) node[] {\pgfuseplotmark{*}};
\end{tikzpicture}
\end{center}
\caption{The graph of the polynomial $p(Z).$}
\label{fig:pol}
\end{figure}
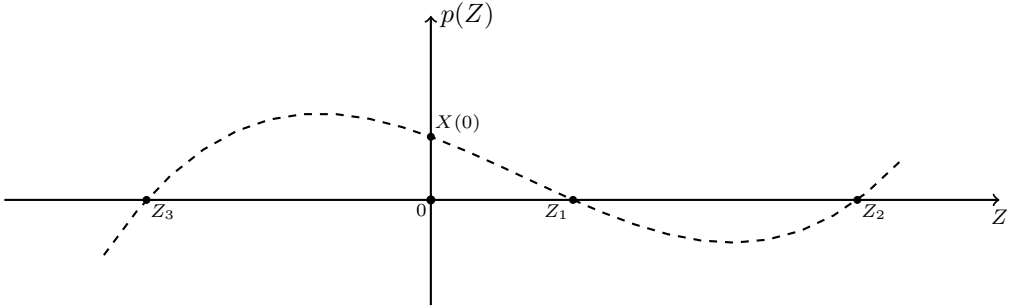

\noindent
For $\, Z\geq 0, \, p(Z)>0 \,$ on the intervals $\, (0,Z_1) \,$ and $ \,(Z_2,\infty). $ Hence, since $ Z(0)=X(0) < Z_1 , $ and $ Z(t) $ is a continuous function, we have by \eqref{dtyest3} that $ Z(t) \leq Z_1 \leq 2 X(0) $ for all $ t. $ This gives \eqref{lemdifineq1}.
\end{proof}

In the proofs below, $ C, C', C_1, C_2, \ldots, $ are numerical constants, independent, in particular, of the parameters $T, K, \eps $ and the function $\rho. $ Subsequent constants depend on the preceding ones at most polynomially.

Now, we recall that by \eqref{defX12} and \eqref{def-x1c},
$X_0(t)=\tfrac12\|\bu^\perp(t)\|^2$ and $ X_1(t) \simeq \|\bu^\perp(t)\|^2
_{1,0}. $ We have
\begin{proposition}\label{prX0est}Assume $\|\bu^\perp(0)\|^2\leq 1/(2C K). $ Then we have
\begin{equation}\label{X0d}
X_0(t)\le 2\, e^{-\s t} \, X_0(0)
\end{equation}
where $ \s := \, 2\,\min (c/(4r_*^2),\, \g\eps), $ with, recall, $ c = 1+\|\rho_x\|^2_{L^\infty}.$
\end{proposition}
\begin{proof} Using that 
$$\frac{c}{4 r_*^2} \, \|u_1^\perp\|^2 \,+\, \gamma \,\|u_2^\perp\|^2\geq \s X_0$$ and
$ \| u^\perp_1\|^4\le 4 X_0^2,\, $ we derive from \eqref{dotx0est}
\begin{align}
\label{X-di} \p_t  X_0 &\le -\frac14\,\|\n u_1^\perp\|^2 \,-\,  \s X_0 \,+\,4\,C_1\,K^2\,X^2_0\,.
\end{align} 
Applying to this inequality Lemma \ref{lemdifineq}, with $\nu =\sigma, \,M_1=4C_1 K^2 $ and $M_2=0, $ yields \eqref{X0d}. 
\end{proof}

\begin{corollary}
Recalling the definitions $X_0\equiv X_0(t):=\frac12 \|\bu^\perp(t)\|^2$ and $\bu^\perp(t):=P^\perp \bu(t)=\bu(t)-\bar \bu(t)$, we arrive at the inequality
\begin{equation}\label{rad-conv-L2}
\|\bu(t)- \bar\bu(t) \| 
\;\le \;C \;e^{-\s t} \; \|\bu_0-\bar \bu_0\|  \qquad (t\ge 0).
\end{equation}
\end{corollary}

To prove a stronger bound, under the additional assumption $ r_*^2/c\leq C_2\gamma/4C_3, $ we define the new, combined Lyapunov function
\begin{equation}\label{definx}
X \,:=\, X_1\;+\;C'\,K^4\, X_0 \;,\qquad \mbox{with} \quad
C'\geq 2C_2/\gamma.
\end{equation}

\noindent
Differentiating $ \, X\,$ and using \eqref{dotx0est} and \eqref{dotx1est}, and absorbing the term $C_1 \,K^2 \, \|u_1^\perp\|^4$ into $C_2 \,K^2\,X_1^2$, we find 
\begin{equation}\label{dotx12est}
\partial_t X
\leq -\, \tilde{X}  
\;+\; C_4\,K^2 \,X^2 
\;+\; C_5\, X^3
\;-\;\frac{1}{8} \, \|\Delta u_1^\perp\|^2\; ,
\end{equation}
\noindent where, with $C' $ as in \eqref{definx},
$$\tilde{X} :=\frac{C'\,K^4}{4} \;\|\nabla u_1^\perp\|^2
\; + \; 
\frac{c\,C'\,}{8\, r_*^2} \; K^4 \; \|u_1^\perp\|^2
\;+\;{\frac12\,\gamma\,C'}\, K^4\, \|u_2^\perp\|^2\,.$$
Now, assuming for simplicity (recall that $ c := 1+\|\rho_x\|^2_{L^\infty}$)
\begin{equation}\label{assmpt}
\max(\eps,\al)\gamma^2\le 2C_2\,,
\end{equation}
and using \eqref{definx}, \eqref{def-x1b} and the definitions of $X_0$ and $ X_1, $ we obtain
\begin{subequations}\label{Xestt}
\begin{eqnarray}
X &\leq& 
\frac12\,\| \nabla u_1^\perp \|^2
\,+\,\frac12 C'K^4\, \|u_1^\perp\|^2
\,+\,\frac12 (\gamma + C' K^4 \eps^{-1})\, \|u_2^\perp\|^2\qquad
\label{Xesti}
\\[0.5em]
&\leq&
\frac{2}{\eps\,\gamma}\,\tilde{X}\;. \label{Xest} \quad\qquad\quad
\end{eqnarray}
\end{subequations}
Dropping the last term in \eqref{dotx12est} and using \eqref{Xest}, we conclude that
\begin{equation}\label{dotx12est2}
\partial_t X \leq - \tfrac{\gamma}{2}\,\eps \, X \;+\;  C_4 \,K^2 \,X^2 
\;+\;C_5 \,X^3 \,,
\end{equation}
where $C_4 $ and $C_5 $ are positive constants which appeared in \eqref{dotx12est} and which are  independent of the parameters $ T, K$ and $ \eps.$

Note that by definition of $C_2 $ in \eqref{e28225}, and for $\gamma\ll 1, $ the second inequality in \eqref{definx} implies the first one and the last two inequalities in \eqref{assmpt}.

\begin{proof}[Proof of Theorem \ref{thm:rad-conv}]\label{proofTm1} Since the Lyapunov function \,$ X(t),\, $ defined in \eqref{definx}, satisfies the differential inequality \eqref{dotx12est2}, by Lemma \ref{lemdifineq} with $\nu=\tfrac{\gamma}{2}\eps, \, M_1=C_4 K^2 $ and $M_2=C_5, $ it obeys estimate \eqref{lemdifineq1}, provided $ X(0) \ll \min\left((C_4 K^2)^{-1}, C_5^{-1}\right)$, which, together with \eqref{def-x1c} and \eqref{definx}, implies \eqref{rad-conv}.
\end{proof}

\subsection{A priori estimates: Proof of (\ref{dotx0est}) and (\ref{dotx1est})}\label{sec:apriori}
In this subsection we prove estimates \eqref{dotx0est} and \eqref{dotx1est} used in subsection \ref{sec28225}. With the function $ X_0(t) $ given in \eqref{defX12}, we have:

\begin{proposition}\label{prop:domain-A}
Let $ T>0, K \geq 1, $ and $\rho$ be a positive, bounded and radially symmetric function of class $C^2$  satisfying $0 < \rho \le r_* $ with $\frac{c}{ r_*^2}\geq 3, $ where, recall, $ c = 1+\|\rho_x\|^2_{L^\infty} \,. $  Then any solution $\bu\equiv \bu(t)=(u_1, u_2)\equiv (u_1(t), u_2(t))$ of \eqref{FHNrho1} in $ S_{T,K}$ obeys the inequality \eqref{dotx0est}, i.e.,
\begin{equation}
\partial_t X_0 \le 
-\,\frac14 \, \|\n u_1^\perp\|^2 \,-\,\frac{c}{4 r_*^2} \, \|u_1^\perp\|^2 \,-\, \gamma \,\|u_2^\perp\|^2
\,+\, C_1 \,K^2 \, \| u^\perp_1\|^4,\label{u-perp-dia} 
\end{equation}
where $C_1$ is a positive constant independent of the parameters $ T, K $ and $\eps $ and the function $\rho.$\end{proposition}

\begin{proof} Using Eq.~\eqref{duperp}, we compute
\begin{align}
\label{derivLyap} \frac12 \p_t \langle  \bu^\perp, \bu^\perp \rangle &= \langle A^\perp \bu^\perp, \bu^\perp \rangle +   \langle \mcN^\perp(\bu), \bu^\perp \rangle.\end{align}
The first term on the r.h.s.\ of \eqref{derivLyap} was estimated in \eqref{Arho-perp-est}.

For the second term on the r.h.s.\ of \eqref{derivLyap}, using definition \eqref{N}, we transform  
\begin{align}
\label{N-perp-ip} \langle \mcN^\perp(\bu), \bu^\perp \rangle = \langle P^\perp_1 h(u_1), u^\perp_1 \rangle \,,
\end{align}
where $ P^\perp_1 $ is the first component of $P^\perp \; (P^\perp= P^\perp_1 \oplus P^\perp_2), \, $ and $ h(u_1)$ is given in \eqref{N}. Next, using that $P^\perp_1 h(\bar u_1)=0$ and $P^\perp_1 h'(\bar u_1) u_1^\perp=h'(\bar u_1) P^\perp_1 u_1^\perp=h'(\bar u_1) u_1^\perp$, we expand 
\begin{align}
\label{g-exp} P^\perp_1 h(u_1)&=P^\perp_1 h(\bar u_1 + u_1^\perp)=h'(\bar u_1) u_1^\perp + h^\perp_1 (u_1),
\\[0.5em]
\label{g1perp}h^\perp_1 (u_1)&:=\frac12 h''(\bar u_1)P^\perp_1  (u^\perp_1)^2+\frac16 h'''(\bar u_1) P^\perp_1  (u^\perp_1)^3.
\end{align}
Combining equations \eqref{N-perp-ip} and \eqref{g-exp}, together with \eqref{g1perp} and the expressions
\begin{equation}\label{gder}
h'(u_1):= -3 u_1^2 +2(\alpha+1)u_1,\;\; h''(u_1):= -6 u_1 +2(\alpha+1), \;\; h'''(u_1):= -6 \,,
\end{equation}
(see \eqref{N}), we arrive at
\begin{eqnarray}
\langle \mcN^\perp(\bu), \bu^\perp \rangle 
& = &
\langle h'(\bar u_1) u_1^\perp, u^\perp_1 \rangle + \langle h^\perp_1(u_1), u^\perp_1 \rangle 
\nonumber\\
& =  &
\langle (-3 \bar u_1^2 +2(\alpha+1)\bar u_1) u_1^\perp, u^\perp_1 \rangle
\nonumber
\\
& &
 \; +
\langle (-3\bar u_1+\alpha+1) (u^\perp_1)^2 - (u^\perp_1)^3, u^\perp_1 \rangle.\label{Nip-exp}
\end{eqnarray}
Using \eqref{Nip-exp}, and using that $ \sup (-3\lambda^2+2(\alpha+1)\lambda)\leq (\alpha+1)^2/3 \, \leq 3/4 $ (recall that $\al\in(0,\tfrac12)$) for the first line in the r.h.s.\ of \eqref{Nip-exp} and dropping the negative term $\langle - (u^\perp_1)^3, u^\perp_1 \rangle$ from the second line, we find
$$\langle \mcN^\perp(\bu), \bu^\perp \rangle \leq \frac34 \|u_1^\perp\|^2\;+\;\frac32 \|u^\perp_1\|_{L^3}^3\;+\;3\int_U|\bar u_1(u_1^\perp)^3|\,d\mu\;.$$
Next, using H\"older estimate 
$ \int_U |\bar{u}_1| |u_1^\perp|^m \, d\mu\leq \|\bar{u}_1\|_{L^\infty}\|u_1^\perp\|_{L^m}^m\, $ and the one-dimensional Sobolev embedding inequality $
 \|\bar{u}_1\|_{L^\infty}  \leq M \|\bar{u}_1\|_{H^1}  $
and recalling that $ \al \le 1, $ we find
\begin{align}\label{Nip-est} 
\langle \mcN^\perp(\bu), \bu^\perp \rangle 
&\le \frac43 \,\|u^\perp_1\|^2 
\; +\; 3 (1+M \|\bar{u}_1\|_{H^1})\,\|u^\perp_1\|_{L^3}^3 \,.
\end{align}
Using \eqref{derivLyap}, \eqref{Arho-perp-est} and \eqref{Nip-est}, we obtain
\begin{align}
\notag \frac12 \p_t \langle  \bu^\perp, \bu^\perp \rangle &\le -\frac12\|\n u_1^\perp\|^2 - \frac{c}{2r_*^2} \|u_1^\perp\|^2 - \gamma \|u_2^\perp\|^2\\
\label{u-perp-di'}&\qquad +\;\frac43 \,\|u^\perp_1\|^2 
\; +\; 3(1+ M \, \|\bar{u}_1\|_{H^1})\,\|u^\perp_1\|_{L^3}^3
\, .\end{align}
Next, we use the Gagliardo-Nirenberg-Sobolev inequality in the dimension $d$:
\begin{align}\label{GNSineq} 
 &\| f\|_{L^m} \le M \|\p^b f\|_{L^2}^{a} \| f\|_{L^2}^{1-a}, \qquad 
 a=  \frac{d}{b}(\frac{1}{2}-\frac{1}{m}),\ m\ge 2, 
\end{align} 
with $b=1, d=2 $ and $ m=3, $ and therefore $ a = \frac13$. Applying this inequality to the last term in \eqref{u-perp-di'} and using 
\begin{subequations}
\begin{eqnarray}
3\; ( M \|\bar u_1\|_{H^1} + 1)\,\|u^\perp_1\|_{L^3}^3 
&\leq &
3\,( M \|\bar u_1\|_{H^1} + 1) \,\|\n u^\perp_1\| \,\| u^\perp_1\|^2
\label{bubperp3a}
\\[0.5em]
&\leq&
\frac{1}{4} \,\|\n u^\perp_1\|^2 
+ C_6\, \big(\|\bar u_1\|_{H^1} + 1\big)^2 \, \| u^\perp_1\|^4 ,\quad
\label{bubperp3b}
\end{eqnarray}
\end{subequations}
where $C_6=(3\max(M,1))^2$, yields the estimate
\begin{align}
\notag \frac12 \, \partial_t \|\bu^\perp\|^2 \le& -\frac14\|\n u_1^\perp\|^2 -\frac{c}{2r_*^2} \|u_1^\perp\|^2-\gamma \|u_2^\perp\|^2
\\[0.5em]
\label{u-perp-di}& 
+ \frac43 \, \| u^\perp_1\|^2 
\,+\, C_6\,\big(\|\bar u_1\|_{H^1} + 1\big)^2  \, \| u^\perp_1\|^4\;.
\end{align}
Using that, by our conditions, $\frac{4}{3}\leq \frac{c}{4r_*^2} $ and recalling that $ \bu \in S_{T,K} $ so that 
\begin{equation}\label{assbu1}
\|\bar u_1(t)\|_{H^1}\le K \; , 
\end{equation}
with $ K \geq 1, $ and recalling the notation $X_0\equiv X_0(t):=\frac12 \|\bu^\perp(t)\|^2, $ we see that Eq.~\eqref{u-perp-di} yields \eqref{u-perp-dia}.
\end{proof}

For the Lyapunov function $ X_1(t) $ given in \eqref{defX12}, we prove the following:

\begin{proposition}\label{propdX1}Let $\rho$ be a positive, bounded and radially symmetric function of class $C^2$  satisfying $0 < \rho \le r_*.$ Then for any  $ T>0, \, K \geq 1, $ any solution $\bu\equiv \bu(t)=(u_1, u_2)\equiv (u_1(t), u_2(t)) \in S_{T,K}$ of \eqref{FHNrho1} obeys the inequality \eqref{dotx1est}, i.e.,
\begin{equation}\label{dotxest}
\begin{split}
\partial_t X_1 
\, \leq& 
\; - \; \frac{1}{8}\;
\|\Delta u_1^\perp\|^2 
\;-\;
\left(\frac{1}{16} \,\bigg(\frac{c}{r_*^2}\right)^2\,-\, C_3 \,K^4\bigg)\; \|u_1^\perp\|^2 
\\
&\, + \, C_2 \big( \|u_2^\perp\|^2
\,+\,K^2\,X_1^2 
\,+ \, X_1^3 \big)
\;,
\end{split}
\end{equation}
where $ C_2 $ and $ C_3 $ are positive constants independent of the parameters $ T, R $ and $ \eps $ and of the function $\rho(x).$
\end{proposition}

\begin{proof}[Proof of Proposition \ref{propdX1}] Omitting the argument $ t $ and using Eq.~\eqref{duperp}, that all the functions involved are real, and using the notation $ \dot{\bu}^\perp = \p_t\bu^\perp , $ we compute
\begin{eqnarray}
\p_tX_1 
&=& - \frac12\langle  \dot{\bu}^\perp, A\bu^\perp \rangle
- \frac12 \langle  \bu^\perp, A\dot{\bu}^\perp \rangle 
\nonumber
\\[0.5em]
&=&
- \frac12\langle  A\bu^\perp, A\bu^\perp \rangle
- \frac12 \langle  \bu^\perp, A^2\bu^\perp \rangle 
 \, + \,\mcV(\bu) \;,
\label{dx1est0}
\end{eqnarray}
\noindent where 
\begin{eqnarray}
\mcV(\bu) &:=& - \frac12 \langle \mcN^\perp(\bu), A\bu^\perp \rangle\;
- \frac12  \, \langle  \bu^\perp, A \mcN^\perp(\bu) \rangle 
\nonumber
\\[0.5em]
&=&
- \langle  \mcN^\perp(\bu) \,,\,\tfrac12 (A+A^*)\bu^\perp \rangle \;,\label{nlt}
\end{eqnarray}
where $ A^* $ is the adjoint of $ A $ in the metric $
\langle \cdot , \cdot \rangle $ \, (see \eqref{inner-rho}). To estimate the nonlinear term $ \mcV, $ we use the Cauchy-Schwarz inequality in \eqref{nlt} to find
\begin{equation}
|\mcV(\bu)|
\,\leq\,
\,  \|\mcN^\perp(\bu)\| \; 
\|\tfrac12\, (A \,+\, A^*) \bu^\perp\|\;.
\label{nlt1}
\end{equation}

Now, using \eqref{A}, we compute \begin{equation}\label{adjA}
\frac12 \big(A+A^*\big) 
\;=\; 
A 
\,+\,
\begin{pmatrix} 0 & 1 \\
-\eps & 0
\end{pmatrix}\;.
\end{equation}

\noindent
Using \eqref{nlt1} and \eqref{adjA}, the triangle inequality, and the relation

\begin{equation}\label{rem}
\Bigg\|{\scriptscriptstyle \begin{pmatrix} 0 & 1 \\
-\eps & 0 
\end{pmatrix}}\bu^\perp\Bigg\|^2
\,=\,
\|u_2^\perp\|^2 \; + \; \eps \; \|u_1^\perp\|^2
\,=\,\eps \, \,\|\bu^\perp\|^2 \; ,
\end{equation}
we obtain furthermore
\begin{eqnarray}
|\mcV(\bu)|
&\leq&
\|\mcN^\perp(\bu)\| \, \big(\|A\bu^\perp\| \, + \,\sqrt{\eps} \, \|\bu^\perp\| \big)
\nonumber
\\[0.5em]
&\leq&
\frac54\,\|\mcN^\perp(\bu)\|^2 \;+\;\frac{1}{4} \; \|A\bu^\perp\|^2 \;+ \;\eps \, \|\bu^\perp\|^2 \; . \label{nltv33}
\end{eqnarray}

\noindent
Using the estimate \eqref{nltv33} in \eqref{dx1est0}, we find
\begin{equation}
\label{dx1est2}
\p_tX_1 \;\leq \;
- \frac14
\|A\bu^\perp\|^2 \;-\; \frac12 \langle  \bu^\perp, A^2\bu^\perp \rangle
\; + \; \frac54\, \|\mcN^\perp(\bu)\|^2 
\;+ \;\eps \, \|\bu^\perp\|^2 \; .
\end{equation}
Now, we estimate the terms in the r.h.s.\ of \eqref{dx1est2}. 
Using \eqref{A}, the expression for the norm corresponding to the inner product \eqref{inner-rho}, i.e.~\eqref{normk0} for $ k=0, $ and the Cauchy-Schwarz inequality, we calculate
\begin{subequations}\label{5e26225}
\begin{flalign}
\|A\bu^\perp\|^2
= \; &
\| (\Delta-\alpha)u^\perp_1-u^\perp_2\|^2 
+ \eps^{-1}\| \eps (u^\perp_1 - \gamma u^\perp_2)\|^2
\nonumber
\\[0.5em]
=\; & \|(\Delta-\alpha)u^\perp_1\|^2
+\| u^\perp_2\|^2 + \eps{\| u^\perp_1\|}^2 +\eps \gamma^2\,\|u^\perp_2\|^2
\nonumber
\\[0.5em]
& -2
\langle (\Delta-\alpha)u^\perp_1, u^\perp_2\rangle
- 2 \eps \gamma\, \langle u^\perp_1, u^\perp_2 \rangle
\nonumber
\\[0.5em]
\geq \;&
\| (\Delta-\alpha)u^\perp_1\|^2
+\| u^\perp_2\|^2
+\eps\| u^\perp_1\|^2
+\eps\gamma^2\| u^\perp_2\|^2
\nonumber
\\[0.5em]
& - 2
\| (\Delta-\alpha)u^\perp_1\|
\| u^\perp_2\|
- 2 \eps \gamma \| u^\perp_1\| \| u^\perp_2\|
\label{auau}
\\[0.5em]
\geq \;&
\| (\Delta-\alpha)u^\perp_1\|^2
+\| u^\perp_2\|^2
+\eps\| u^\perp_1\|^2
+\eps\gamma^2\| u^\perp_2\|^2
\nonumber
\\[0.5em]
& - 2 \| (\Delta-\alpha)u^\perp_1\|^2
- \,\frac12 \, \| u^\perp_2\|^2
- \, 2 \eps \, \| u^\perp_1\|^2
- \, \frac{\eps\gamma^2}{2} \,\| u^\perp_2\|^2\;.\label{auaub}
\end{flalign}
\end{subequations}
Collecting similar terms, we conclude
\begin{equation}
\|A\bu^\perp\|^2
\;\geq \;  - \, \|(\Delta-\alpha)u^\perp_1\|^2
\;-\; \eps \, \|u^\perp_1\|^2 
\; + \;
\big(\frac12 + \frac{\eps\gamma^2}{2}\big) \, \| u^\perp_2\|^2 \;.
\label{L2A}
\end{equation}
For the second term in the r.h.s.\ of \eqref{dx1est2}, we claim that
\begin{equation}\label{ua2u}
\langle  \bu^\perp, A^2\bu^\perp \rangle
\;=\;
\|(\Delta-\alpha)u^\perp_1\|^2
\;-\; \eps \, \|u^\perp_1\|^2
\; - \;
(1 - \eps\gamma^2)
\|u_2^\perp\|^2\;.
\end{equation}
Indeed, using the relation
\begin{equation}\label{asquare-eq1}
A^2 
 = 
\begin{pmatrix} \Delta -\al & -1 \\
\eps & - \eps\gamma 
\end{pmatrix}^2
 =
\begin{pmatrix} 
(\Delta -\al)^2 -\eps & -(\Delta -\al)+ \eps \gamma  
\\
\eps (\Delta - \al) - \eps^2\gamma & - \eps +(\eps\gamma)^2 
\end{pmatrix}
\end{equation}

\noindent  and definition \eqref{inner-rho}, we obtain

\begin{flalign}
\langle  \bu^\perp, A^2\bu^\perp \rangle
= \; &
\langle u_1^\perp \,,\, [(\Delta -\al)^2 -\eps] u_1^\perp - [(\Delta -\al)-  \eps \gamma] u_2^\perp \rangle
\nonumber
\\[0.5em]
&\; + 
\frac1\eps \langle u_2^\perp \,,\, [\eps (\Delta - \al) - \eps^2\gamma] u_1^\perp  -\eps (1 -\eps \gamma^2 ) u_2^\perp\rangle
\nonumber
\\[0.5em]
= \; & 
\langle u_1^\perp  \, , \, (\Delta-\al)^2 u_1^\perp\rangle
-\eps \|u_1^\perp\|^2
-
\langle u_1^\perp  \, , \, (\Delta -\al) u_2^\perp\rangle
\nonumber
\\[0.5em]
& 
+\eps\gamma
\langle u_1^\perp  \, , \, u_2^\perp\rangle
+
\langle u_2^\perp  \, , \, (\Delta -\al) u_1^\perp\rangle
\nonumber
\\[0.5em]
& \;
-\eps\gamma
\langle u_2^\perp  \, , \, u_1^\perp\rangle
- (1 - \eps\gamma^2)
\|u_2^\perp\|^2\;.
\label{upa2up-5}
\qquad\qquad\qquad
\end{flalign}
Now, collecting similar terms and using that $ \Delta $ is self-adjoint on $ L^2 $ and therefore
\begin{equation}\label{u1d2u1}
\langle u_1^\perp  \, , \, (\Delta-\al)^2 u_1^\perp\rangle
\;=\;\|(\Delta-\al) u_1^\perp\|^2 \;,
\end{equation}
we arrive at \eqref{ua2u}. Adding inequalities \eqref{L2A} and \eqref{ua2u} we arrive at
\begin{equation}
\frac12\,\|A\bu^\perp\|^2
\,+\,
\langle  \bu^\perp, A^2\bu^\perp \rangle
\,\geq\,
\frac12 \,  \|(\Delta-\alpha)u^\perp_1\|^2
\,-\, \frac{3\eps}{2}\, \|u_1^\perp\|^2
\,- \, \big(\frac34 \,-\, \frac{ 5 \,\eps\, \gamma^2}{4}\big) \,
\|u_2^\perp\|^2 \,.\qquad
\label{17524}
\end{equation}
Next, we estimate $ \|(\Delta-\alpha)u^\perp_1\|. $ Using that $ \Delta $ is self-adjoint on $ L^2, $ we find
\begin{equation}\label{17524b}
\|(\Delta-\alpha)u^\perp_1\|^2
\;=\;
\|\Delta u^\perp_1\|^2
\; - \; 
2\, \al \,
\langle u^\perp_1 , \Delta u^\perp_1\rangle
\; + \;
\al^2\,
\| u^\perp_1\|^2 \;.
\end{equation}
Writing $\, \|\Delta u^\perp_1\|^2=\langle\sqrt{-\Delta}u^\perp_1\,,\,(-\Delta)\sqrt{-\Delta}u^\perp_1\rangle$ and using \eqref{Lapl-est} twice, we obtain
\begin{equation}\label{ulugequ}
\|\Delta u^\perp_1\|^2
\;\geq \; 
\left(\frac{c}{r_*^2}\right)^2\;
\|u^\perp_1\|^2\;,
\end{equation}
where, recall, $ c=1+\|\partial_x\rho\|_{L^\infty}. $ This, together with \eqref{17524b}, gives
\begin{equation}\label{Laplau1}
\|(\Delta-\al)u^\perp_1\|^2 
\;\geq\;
\frac12\,\|\Delta u^\perp_1\|^2
\;+\;2\al \,\|\nabla u_1^\perp\|^2
\;+\;
\big(\frac{c^2}{2 r_*^4} \,+\, \al^2\big)\; \|u^\perp_1\|^2 
\; . 
\end{equation}
Inequalities \eqref{17524} and \eqref{Laplau1} and the assumption $\frac{c^2}{4 r_*^4} \,+\, \al^2\geq 7 \eps$ imply
\begin{equation}
\begin{split}
\frac12\,\|A\bu^\perp\|^2
\,+\,
\langle  \bu^\perp, A^2\bu^\perp \rangle
\, \geq &\;
\frac{1}{4} \,  \|\Delta u^\perp_1\|^2
\;+\;\al \,\|\nabla u_1^\perp\|^2
\\[0.5em]
& + \, \left( \frac{c^2}{8 \, r_*^4}\,+\,2\eps\right) \, \|u_1^\perp\|^2
\,- \, 2 \, C_7\,
\|u_2^\perp\|^2 \,,\qquad
\label{17524c}
\end{split}
\end{equation}
where
\begin{equation}\label{17524d}
C_7 =  \frac{3}{8} \,- \, \frac{5}{8}\, \; \eps \, \gamma^2 \;.
\end{equation} 
Using inequality \eqref{17524c} in \eqref{dx1est2} and writing $\, \|\bu^\perp\|^2 \,=\, \|u_1^\perp\|^2+\eps^{-1}\|u_2^\perp\|^2 \, ,\,$  we find the estimate
\begin{equation}
\partial_t X_1 
\, \leq 
\,-\frac{1}{8} \;
\|\Delta u_1^\perp\|^2
\;-\;
\frac{1}{16}\left(\frac{c}{r_*^2}\right)^2 \, \|u_1^\perp\|^2
\; + \;
C_8 \,
\|u_2^\perp\|^2
\; + \; \frac54\;\|\mcN^\perp(\bu)\|^2 \;,
\label{dotx1est1}
\end{equation}
with
\begin{equation}\label{c2p} C_8=C_7+1, \end{equation} 
where \,$ C_7 $ is defined in \eqref{17524d}. Note that, with \eqref{17524d}, \eqref{c2p} and $\gamma\ll 1, $ the last condition in \eqref{definx} implies the middle condition in \eqref{definx} and the last two conditions in \eqref{assmpt}.

Finally, we estimate the norm $ \|\mcN^\perp(\bu)\| $ of the nonlinearity $\mcN^\perp(\bu)$ entering \eqref{dotx1est1}.

\begin{lemma}\label{lemmaN}
\begin{equation}\label{np2a}
\|\mcN^\perp(\bu)\|^2 \leq C_9\, K^4\, \|u_1^\perp\|^2 \, +\, C_{10} \sum_{ r=4,6}K^{6-r} \|u^\perp_1\|_{H^1}^r \;.
\end{equation}
\end{lemma}

\noindent
We prove this lemma at the end of this section. Now, we proceed with the proof of Proposition~\ref{propdX1}. Using Eqs \eqref{def-x1b} and \eqref{np2a} we arrive at the inequality
\begin{equation}\label{np3}
\|\mcN^\perp(\bu)\|^2 \;\leq \;C_{11} \, \big(K^4 \, \|u_1^\perp\|^2 
\, + \, K^2\, X_1^2 \, +\, X_1^3 \big)\,,
\end{equation}
for some constant $ \,C_{11} > 0  \, $ independent of $K.$ Inserting \eqref{np3} into \eqref{dotx1est1} and using \eqref{ulugequ}, we find
\begin{subequations}\label{dX1estsub}
\begin{flalign}
\dot{X}_1 
 \leq &
\,-\frac{1}{8} \;
\|\Delta u_1^\perp\|^2 
\;-\;
\bigg(\frac{1}{16}\bigg(\frac{c}{r_*^2}\bigg)^2\,-\,C_3\, K^4\bigg) \, \|u_1^\perp\|^2
\label{dX1est3}
\\[0.5em]
& 
\, + \, C_2 \big( \|u_2^\perp\|^2
\,+\,K^2\,X_1^2 
\,+ \, X_1^3 \big) \;.
\label{dX1est4}
\end{flalign}
\end{subequations}
This gives \eqref{dotxest}, completing the proof of Proposition \ref{propdX1}.
\end{proof}

\begin{remark}\label{3r27225}
Equation \eqref{ua2u} shows that the Lyapunov function $ X_2(t)$
$:=\frac12 $ $\langle  \bu^\perp(t), (-A)^2\bu^\perp(t) \rangle $ (cf.~\eqref{defX12}) is positive modulo zero order terms.
\end{remark}

\begin{proof}[Proof of Lemma \ref{lemmaN}]
Using \eqref{N} and \eqref{g-exp}-\eqref{g1perp}, we obtain
\begin{flalign}
\|\mcN^\perp(\bu)\|
=
\|P_1^\perp h(u_1)\|
\leq \;&
\|h'(\bar u_1)\|_{L^\infty}\|u_1^\perp\|
\,+\,
\tfrac12\,
\|h''(\bar u_1)\|_{L^\infty} \|P_1^\perp(u_1^\perp)^2\|
\nonumber
\\[0.5em]
&
+\,
\tfrac16\|h'''(\bar u_1)\|_{L^\infty}
\|P_1^\perp(u_1^\perp)^3\|\,.\label{nperp1}
\end{flalign}
We estimate each term on the r.h.s.\ of \eqref{nperp1}. Using \eqref{gder} and estimating $ \, \|h^{(k)}(\bar u_1)\|_{L^\infty}\, $ in terms of $ \, \|\bar u_1\|_{L^\infty} , \, $ and then applying Sobolev embedding inequality for $ \bar u_1, $ we find
\begin{eqnarray}
 \|h'(\bar u_1)\|^2_{L^\infty} 
& \leq &  
C_{12} \big(\|\bar u_1\|^4_{L^\infty} + \|\bar u_1\|_{L^\infty}^3 +\|\bar u_1\|^2_{L^\infty}\big)
\nonumber
\\[0.5em]
& \leq &
C_{13} \big(\|\bar u_1\|^4_{H^1} + \|\bar u_1\|_{H^1}^3 +\|\bar u_1\|^2_{H^1}\big) \; ,
\label{gdinfty}
\\[1em]
\|h''(\bar u_1) \|_{L^\infty}^2 
& \leq &
C_{14} \, 
\big(
\|\bar u_1\|_{L^\infty}^2
+
\|\bar u_1\|_{L^\infty}
+ 1
\big)
\nonumber
\\[0.5em]
&\leq&
C_{15}  \,
\big(
\|\bar u_1\|_{H^1}^2
+
\|\bar u_1\|_{H^1}
+1
\big)  \;,
\label{g2dinfty}
\end{eqnarray}
\noindent
and $ \, \|h'''(\bar u_1) \|_{L^\infty} = 6. $

Next, for the term $ \|P_1^\perp(u_1^\perp)^2\| $ in the second term of the r.h.s.\ of \eqref{nperp1}, we use that $ \|P_1^\perp\| \leq 1 $ to obtain
\begin{equation}\label{p1p2}
\|P_1^\perp(u_1^\perp)^2\|^2
\,\leq\, 
\|(u_1^\perp)^2\|^2
\,=\,
\|u^\perp_1\|_{L^4}^4\,.
\end{equation}
Similarly, for the term $ \|P_1^\perp(u_1^\perp)^3\| $ in \eqref{nperp1}, we estimate $ \|P_1^{\perp}((u^\perp_1)^3)\|_{L^2}^2
\leq \|(u^\perp_1)^3\|_{L^2}^2 $ 
$=$
$ \|u_1^\perp\|_{L^6}^6 $ to find 
\begin{equation}\label{p1p3}
\|P_1^\perp(u_1^\perp)^3\|^2
\; \leq\;
 \|u_1^\perp\|_{L^6}^6\;.
\end{equation}
\noindent 
Inserting \eqref{gdinfty}, \eqref{g2dinfty}, \eqref{p1p2}, \eqref{p1p3} into \eqref{nperp1}, we obtain
\begin{equation}\label{np4}
\|\mcN^\perp(\bu)\|^2 \leq 
C_{16} \, \sum_{k=2,3,4}\,\|\bar u_1\|_{H^1}^k \, \|u_1^\perp\|^2 \, + \, C_{17}\sum_{k+r=4,5,6, 
\atop r=4,6}\,\|\bar u_1\|_{H^1}^k \; \|u^\perp_1\|_{L^r}^r 
\end{equation}
and using that  \;$ \|u_1^\perp\|_{L^r} < M\, \|u_1^\perp\|_{H^1} \,, \;\, r=4,6, $ we find
\begin{equation}\label{np2}
\|\mcN^\perp(\bu)\|^2 \leq C_{18}  \, \sum_{\scriptscriptstyle k=2,3,4}\,\|\bar u_1\|_{H^1}^k \, \|u_1^\perp\|^2_{L^2} \, +\, C_{19}\sum_{k+r=4,5,6, 
\atop r=4,6}\,\|\bar u_1\|_{H^1}^k \|u^\perp_1\|_{H^1}^r \;.
\end{equation}
Recalling the assumption $ \|\bar{u}_1\|_{H^1}\leq K ,\;K\gtrsim 1, $ and using \eqref{np2}, we arrive at \eqref{np2a}.
\end{proof}

\section{Effective equation: Proof of Theorem \ref{thm:average}}\label{seqefeq}

First, we derive a convenient equation for $ \bar{\bm{u}}. $ 
Applying the projection $ \,P\,$ (see \eqref{proj}) to \eqref{rho-eq-canon} and using decomposition \eqref{decomp},  we find
\begin{equation}\label{deqv}
\p_t \bar{\bm{u}}\; = \;
\Ar \bar{\bm{u}} \;+ \;P \mcN(\bar{\bu} + \bu^\perp)\; ,
\end{equation} 
where, recall, $\, \mcN(\bu) \, $ is given by \eqref{N} and $ \, \Ar\,$ is defined as (cf.~\eqref{A})
\begin{equation}\label{Arad} 
\Ar  
:=
\begin{pmatrix} \Dr -\alpha & -1 \\
\eps & - \eps\gamma 
\end{pmatrix}\,.
\end{equation}
Here we used the relation $ \, PA=\, \Ar P \, $ which follows from the easy verifiable equations $\, P A = A P \, $ and $ \, A P = \Ar P .$

Furthermore, we rewrite Eq.~\eqref{eq:radlapl} as
\begin{equation}\label{eqbw}
\p_t\bm{w} = \Ar \bm{w}\,+\,\mcN(\bm{w}) \;,
\end{equation}
with $ \, \Ar \, $ and $ \, \mcN(\bm{w}) \, $ given in \eqref{Arad} and \eqref{N}.

Next, we let $\,\bm{v}(t) = \bar\bu(t)-\bm{w}(t) , \,$  and define the Lyapunov functional
\begin{equation}\label{X1diag}
Y_1(t) \, = \, -\, \frac12\, 
\big\langle A_d\bm{v}(t)\,,\,\bm{v}(t)\big\rangle \;,
\end{equation}
where $ \,A_d \, $ is the diagonal part of the operator $ \Ar$ :
\begin{equation}\label{Ad}
A_d  
\;:=\;
\begin{pmatrix} \Dr -\alpha & 0 \\
0 & - \eps\gamma 
\end{pmatrix}
\;.
\end{equation}
Integrating by parts in \eqref{X1diag}, we obtain readily
\begin{equation}\label{X1dest}
0\leq Y_1(t) \simeq \|\bm{v}(t)\|_{1,0}^2 \;.
\end{equation}
Notice the difference between \eqref{X1diag} and \eqref{defX12} for $ k=1 $ (or \eqref{def-x1}).

\begin{proposition}\label{vdifineq}
Let $\bu(t)$ and $\bm{w}(t) $ satisfy Eqs.~\eqref{FHNrho1} and \eqref{eq:radlapl}, respectively. Then, the function $ Y_1(t) $ defined in \eqref{X1diag}, satisfies the differential inequality
\begin{eqnarray}\label{vdifineq0}
\p_t Y_1(t) &\leq&\, C \, \left( Y_1(t) \,+\,  Y_1^3(t) 
\, + \, 
\big\|u_1^\perp(t)\big\|_{H^1}^{4}
\,+\,
\big\|u_1^\perp(t)\big\|_{H^1}^{8}
\right)
\,, 
\end{eqnarray}
for some constant $ \, C>0.$
\end{proposition}

\begin{proof}Various constants appearing in this proof are numerical constants due mainly to Sobolev embedding inequalities, i.e.\ polynomials in the Sobolev constant $M. $ They are independent of the constants in Section \ref{sec:rad-collapse}. Moreover, subsequent constants are at most polynomials in the preceding ones.

Expanding $ \, \mcN(\bar{\bu} +\bu^\perp) \, $ in $ \, \bu^\perp \, $ and using that $ \, P \mcN'(\bar{\bu})\bu^\perp = \mcN'(\bar{\bu})P \bu^\perp = 0 , \, $ we obtain
\begin{equation}\label{PNu}
P \mcN(\bar{\bu} + \bu^\perp) \,=\, \mcN(\bar{\bu})\,+\,\mcR(\bar{\bu}, \bu^\perp)\;,
\end{equation}
with
\begin{equation}\label{Ravuuperp}
\mcR(\bar{\bu}, \bu^\perp) 
\;:=\;
\frac{1}{2!} \, \mcN''(\bar{\bu}) \, P(\bu^\perp \bu^\perp)\, + \, \frac{1}{3!} \, \mcN'''(\bar{\bu})\, P(\bu^\perp \bu^\perp \bu^\perp)
\,=\,\mathcal O(\|\bu^\perp\|^2) \;.
\end{equation}
Eqs.~\eqref{deqv}, \eqref{eqbw} and \eqref{PNu}  imply that the vector-function $ \,\bm{v}(t)  \, $ satisfies the equation 
\begin{equation}\label{eqptv}
\p_t\bm{v}(t)\,=\, \Ar \bm{v}(t)\,+\,\mcN(\bar{\bu}) \, +\, \mcR(\bar{\bu}, \bu^\perp) \, - \, \mcN(\bm{w}) \,,
\end{equation}
with initial condition $ \, \bm{v}(0)=0. $ We use \eqref{eqptv} to derive a differential inequality for the Lyapunov functional $ \, Y_1(t) \,$ defined in \eqref{X1diag}. We use  \eqref{X1diag}, \eqref{eqptv} and the fact that the operator $A_d$ is symmetric to obtain
\begin{eqnarray}
\p_t Y_1 &=& 
- \, \frac12 \,
\big\langle A_d\dot{\bm{v}} \, , \, \bm{v}\big\rangle
\;-\;\frac12\,
\big\langle A_d\bm{v} \, , \, \dot{\bm{v}}\big\rangle
\nonumber
\\[0.5em]
&=&
-\,\big\langle A_d\bm{v} \, , \, \dot{\bm{v}}\big\rangle
\nonumber
\\[0.5em]
&=&
- \, \big\langle A_d\bm{v} \, , \, \Ar\bm{v}\big\rangle
\; - \;
\big\langle 
A_d\bm{v}
\; , \;
\mcM(\bar{\bu})
\big\rangle \;,
\label{ptx1est1}
\end{eqnarray}
\noindent where, recall, $ \dot{\bm{v}}=\partial_t\bm{v} $ and
\begin{equation}\label{Mbaru}
\mcM(\bar{\bu})\;:=\; \mcN(\bar{\bu}) \, +\, \mcR(\bar{\bu}, \bu^\perp) \, - \, \mcN(\bm{w}) \;.
\end{equation}
We estimate the terms on the r.h.s.\ of \eqref{ptx1est1}. For the first term in the r.h.s.\ of \eqref{ptx1est1}, we recall the definition \eqref{Arad} of $ \, \Ar, \, $ and we write
\begin{equation}\label{Asplit} 
\Ar  
\,=\,
\begin{pmatrix} \Dr -\alpha & -1 \\
\eps & - \eps\gamma 
\end{pmatrix}
\;=\;
A_d\;+\;V \; ,
\end{equation}
where, recall,
\begin{equation}\label{AdV}
A_d  
\;:=\;
\begin{pmatrix} \Dr -\alpha & 0 \\
0 & - \eps\gamma 
\end{pmatrix}
\qquad\mbox{and}\qquad
V\;:=\;
\begin{pmatrix} 0 & -1 \\
\eps & 0 
\end{pmatrix}
\;.
\end{equation}
We use \eqref{Asplit}-\eqref{AdV} and apply Cauchy-Schwarz and Young inequalities to get
\begin{eqnarray}
-\,
\big\langle A_d\bm{v} \, , \, \Ar\,\bm{v}\big\rangle
&=&
-\, 
\big\langle A_d\bm{v} \, , \, (A_d+V)\bm{v}\big\rangle
\nonumber
\\[0.5em]
& = &
-\, \|A_d\bm{v}\|^2
\;-\;
\big\langle A_d\,\bm{v} \, , \, V\bm{v}\big\rangle
\nonumber
\\[0.5em]
& \leq &
-\, \|A_d\bm{v}\|^2
\;+\;
\| A_d\bm{v}\| \, \|V\bm{v}\|
\nonumber
\\[0.5em]
& \leq &
-\,\frac12\, \|A_d\bm{v}\|^2
\;+\;
\frac12 \, \|V\bm{v}\|^2
\nonumber
\\[0.5em]
& = &
-\,\frac12\, \|A_d\bm{v}\|^2
\;+\;
\frac{\epsilon}{2} \, \|\bm{v}\|^2\;.
\label{ptx1est2}
\end{eqnarray}
Using \eqref{ptx1est2} in \eqref{ptx1est1}, we obtain
\begin{eqnarray}
\p_t Y_1 
&\leq & 
-\,\frac12\, \|A_d\bm{v}\|^2
\;+\;
\frac{\epsilon}{2} \, \|\bm{v}\|^2
\; - \,
\big\langle 
A_d\bm{v}
\; , \;
\mcM(\bar{\bu})
\big\rangle\;,
\label{ptx1est3}
\end{eqnarray}
where $ \,\mcM(\bar{\bu})\, $ is defined in \eqref{Mbaru}. To estimate the last term in \eqref{ptx1est3}, we use that, by the definition, $ \,\mcM(\bar{\bu}) = \big(\mcM_1(\bar{u}_1) , 0\big),\, $ and therefore,
\begin{eqnarray}
\big|\big\langle 
A_d\bm{v}
\; , \;
\mcM(\bar{\bu})
\big\rangle\big|&=&
\big|\big\langle 
(\Dr-\al)v_1
\; , \;
\mcM_1(\bar{u}_1)
\big\rangle\big|
\nonumber
\\[0.5em]
&\leq&
\big\|(\Dr-\al)v_1\big\|\;\big\|\mcM_1(\bar{u}_1)\big\|
\nonumber
\\[0.5em]
&\leq&
\frac14 \, \big\|(\Dr-\al)v_1\big\|^2
\, + \, \big\|\mcM_1(\bar{u}_1)\big\|^2\;.\label{AdvMbar}
\end{eqnarray}
Next, writing $ \, \bar{\bu} = \bm{w}+\bm{v} ,\, $ we expand $ \, \mcN(\bar{\bu}) = N(\bm{w}+\bm{v}) \, $ in $ \, \bm{v}\, $ to obtain
\begin{equation}\label{NbaruNw}
\mcN(\bar{\bu}) \;- \; \mcN(\bm{w}) 
\;=\;
\mcN'(\bm{w}) \, \bm{v}
\;+\;
\frac{1}{2!} \, \mcN''(\bm{w})\,\bm{v}\,\bm{v}
\;+\;
\frac{1}{3!}\,\mcN'''(\bm{w})\,\bm{v} \,\bm{v} \,\bm{v}
\end{equation}
and
\begin{eqnarray}
\mcR(\bar{\bu}, \bu^\perp) 
&=&
\mcR(\bm{w}, \bu^\perp)
\;+\;
\mcR'(\bm{w}, \bu^\perp) \, \bm{v}
\nonumber
\\[0.5em]
& & \;+\;
\frac{1}{2!} \,\mcR''(\bm{w}, \bu^\perp)\,\bm{v}\,\bm{v}
\;+\;
\frac{1}{3!}\,\mcR'''(\bm{w}, \bu^\perp)\,\bm{v} \,\bm{v} \,\bm{v}
\label{Rbaruperp}
\end{eqnarray}
where $ \, \mcN'', \mcR'' \, $ and $ \, \mcN''', \mcR''' \, $ are bilinear and trilinear maps (tensors), respectively, and the derivatives of $ \, \mcR \, $ are taken with respect to $\, \bm{w}$.

Furthermore, by definition \eqref{Ravuuperp} and the fact that $\mcN(\bu)=(\mcN_1(u_1),0), $ with $\mcN_1(u_1) $ being a third order polynomial in $u_1, $ we have
\begin{equation}\label{Rder}
\mcR'(\bm{w}, \bu^\perp) 
\;=\;
\frac{1}{2!} \, \mcN'''(\bm{w})\,P(\bu^\perp\bu^\perp)\,,
\quad \mcR''(\bm{w}, \bu^\perp) = 0\,, \quad \mcR'''(\bm{w}, \bu^\perp)=0\;.
\end{equation}
Eqs \eqref{Rbaruperp} and \eqref{Rder}, together with Eq.~\eqref{Ravuuperp} and the fact that only the first components in $\mcR$ and $\mcN$ do not vanish, yield
\begin{equation}\label{Ravuwperp}
\begin{split}
\mcR_1(\bar{\bu}, \bu^\perp) 
\;:=& \;
\frac{1}{2!} \, \mcN_1''(w_1) \, P_1((u_1^\perp)^2)\, + \, \frac{1}{3!} \, \mcN_1'''(w_1)\, P_1((u_1^\perp)^3) \\& + \, \frac{1}{2!} \, \mcN_1'''(w_1)\,P_1((u_1^\perp)^2)\,v_1\;.
\end{split}
\end{equation}
Next, we use \eqref{NbaruNw} and \eqref{Ravuwperp} in \eqref{Mbaru} to compute
\allowdisplaybreaks
\begin{eqnarray}
\mcM_1(\bar{u}_1)
&=&
\mcN_1'(w_1) \, v_1
\;+\;
\frac{1}{2!} \, \mcN_1''(w_1)\,v_1^2
\;+\;
\frac{1}{3!}\,\mcN_1'''(w_1)\,v_1^3
\nonumber
\\[0.5em]
& &
+ \;
\frac{1}{2!} \, \mcN_1''(w_1) \, P_1((u_1^\perp)^2)\, + \, \frac{1}{3!} \, \mcN_1'''(w_1)\, P_1((u_1^\perp)^3)
\nonumber
\\[0.5em]
& & + \;
\frac{1}{2!} \, \mcN_1'''(w_1)\,P_1((u_1^\perp)^2)\,v_1
\;=\; B'\,+\,B''\,+\,B'''\;,
\label{ptx1est4}
\end{eqnarray}
where $B', B''$ and $ B''' $ stand for the first, second and third lines on the r.h.s.\ of \eqref{ptx1est4}.

Here, unlike \eqref{Rbaruperp} and \eqref{Rder}, $ \mcN_1'(w_1), $ etc, are usual derivatives of scalar functions.

Now, we estimate the terms $ \,P_1\big((u_1^\perp)^2\big)\, $ and $ \, P_1\big((u_1^\perp)^3\big)\, $ in \eqref{ptx1est4}. We use definition \eqref{proj} of projection $ \, P, \, $ and apply H\"older inequality, to obtain, for $ \, k=2,3,$
\begin{equation}
\left|P_1\big((u_1^\perp)^k\big)\right|
\;=\;
\frac{1}{2\pi}  \, \left|
\int_0^{2\pi} (u_1^\perp)^k \, d\theta\right|
\;\leq\;
\frac{1}{\sqrt{2\pi}}\,\left(\int_0^{2\pi}(u_1^\perp)^{2k} \, d\theta\right)^{1/2}\;.
\label{puperp21}
\end{equation}
Recall the definition of $ L^2_{rad} $ in \eqref{L2rhorad} and let $L^p \equiv L^p(\cS_{\rho}, \mathbb R)$ (cf.~\eqref{vecL2-rho}). Taking $L^2_{rad}$-norm in \eqref{puperp21} and applying two-dimensional Sobolev embedding inequality, we obtain
\begin{eqnarray}
\Big\|P_1\big((u_1^\perp)^k\big)\Big\|_{L^2_{\rm rad}}
&\leq& 
\frac{1}{\sqrt{2\pi}}\,
\bigg\|\bigg(\int_0^{2\pi}(u_1^\perp)^{2k} \, d\theta\bigg)^{1/2}\bigg\|_{L^2_{ \rm rad}}
\nonumber
\\[0.5em]
&=&
\frac{1}{\sqrt{2\pi}}\,
\bigg(
\int_0^\infty
\bigg|\,
\bigg(
\int_0^{2\pi}(u_1^\perp)^{2k} \, d\theta
\bigg)^{1/2}
\,\bigg|^2
\,\sqrt{g}\,dx
\bigg)^{1/2}
\nonumber
\\[0.5em]
&=&
\frac{1}{\sqrt{2\pi}}\,
\left(
\int_0^\infty
\int_0^{2\pi}(u_1^\perp)^{2k} \,\sqrt{g}\; d\theta
\,dx
\right)^{1/2}
\nonumber
\\[0.5em]
&=&
\tfrac{1}{\sqrt{2\pi}}\;
\big\|u_1^\perp\big\|_{_{L^{2k}}}^k
\nonumber
\\[0.5em]
&\leq&
\tfrac{M^k}{\sqrt{2\pi}} \,\big\|u_1^\perp\big\|_{H^1}^k\;.
\label{puperp22}
\end{eqnarray}

\noindent
Now, we estimate the terms in the r.h.s.\ of \eqref{ptx1est4}. For the first term in the r.h.s.\ of \eqref{ptx1est4}, we use \eqref{N}, \eqref{gder}, and triangle,  H\"older and Sobolev inequalities and the condition $\al< \tfrac12$, to obtain
\begin{flalign}
\big\|B'\big\| &=
\big\|
(3w_1^2-2(\al+1))\,v_1
\;+\;(3w_1-\alpha-1)
\, v_1^2 \;+\; v_1^3 \big\|
\nonumber
\\[0.5em]
&\leq
3 \big(\|w_1\|_\infty^2 + 1\big)
\,\big\|v_1\big\|
\;+\;
3 \big(\|w_1\|_\infty   + 1\big) \, \big\|v_1\big\|^2_{L^4}
\;+\;
\big\|v_1\big\|^3_{L^6}
\nonumber
\\[0.5em]
&\leq
3 \big(\|w_1\|_\infty^2 + 1\big)
\,\big\|v_1\big\|
\;+\;
3 M^2\, \big(\|w_1\|_\infty   + 1\big) \, \big\|v_1\big\|^2_{H^1}
\;+\;
M^3\,\big\|v_1\big\|^3_{H^1}
\;.
\label{ptx1est5}
\end{flalign}

For the second term in the r.h.s.\ of \eqref{ptx1est4}, we use \eqref{N}, \eqref{gder}, \eqref{puperp22}, and triangle and H\"older inequalities, to obtain

\begin{eqnarray}
\big\|B''\big\| &=&
\big\| (3 w_1 - \al-1) \, P_1\big((u_1^\perp)^2\big)
\;+\;
P_1\big((u_1^\perp)^3\big) \big\|
\nonumber
\\[0.5em]
&\leq&
3 \big(\|w_1\|_\infty +  1\big) \;\left\|P_1\big((u_1^\perp)^2\big)\right\|
\;+\;
\left\|P_1\big((u_1^\perp)^3\big)\right\|
\nonumber
\\[0.5em]
&\leq& \tfrac{3 M^2}{\sqrt{2\pi}} \, \big(\|w_1\|_\infty + 1\big) \;\big\|u_1^\perp\big\|_{H^1}^2
\;+\;
\tfrac{M^3}{\sqrt{2\pi}} \, \big\|u_1^\perp\big\|_{H^1}^3
\;.
\label{ptx1est6}
\end{eqnarray}
For the last term in the r.h.s.~of \eqref{ptx1est4}, we use \eqref{N}, \eqref{gder}, \eqref{puperp21}, H\"older, Sobolev and triangle inequalities to obtain
\begin{eqnarray}
\big\|B'''\big\|
&=&
3\; \big\|P_1\big((u_1^\perp)^2\big)
\, v_1\big\|
\nonumber
\\[0.5em]
&\leq &
\frac{3}{(2\pi)^{3/4}} \, \big\|u_1^\perp\big\|_{L^8}^2 \; \big\|v_1\big\|_{L^4}
\label{e23225}
\\[0.5em]
& \leq &
\frac{3M^3}{(2\pi)^{3/4}} \, \big(\big\|u_1^\perp\big\|_{H^1}^4 \,+\, \big\|v_1\big\|_{H^1}^2\big)\;.
\label{ptx1est7p}
\end{eqnarray}
Therefore, by \eqref{ptx1est4} and \eqref{ptx1est5}-\eqref{ptx1est7p}, we obtain
\begin{flalign}\label{e24225}
\big\|\mcM_1(\bar{u}_1)\big\|
\leq 
\big\|B'\big\|+\big\|B''\big\|+\big\|B'''\big\|
 \leq \;
C_1\,
\Big(&
\,\big\|v_1\big\|^2_{H^1}
\,+\,
\big\|v_1\big\|^3_{H^1}
\,+\,
\big\|v_1\big\|
\nonumber
\\[0.8em]
&+\,
\big\|u_1^\perp\big\|_{H^1}^2
\,+\,
\big\|u_1^\perp\big\|_{H^1}^4
\Big)\;.
\end{flalign}
Using \eqref{AdvMbar} and $\;\|A_d\bm{v}\|^2=\|(\Dr -\alpha)v_1\|^2+\eps \gamma^2\|v_2\|^2 \; $ in \eqref{ptx1est3}, we conclude that
\begin{equation}
\p_t Y_1 \leq 
-\,\frac14\,  \big\|(\Dr -\alpha)v_1\big\|^2
\,-\,
\frac{\eps \gamma^2}{2} \, \|v_2\|^2 
\,+\,
\frac{\epsilon}{2} \, \big\|\bm{v}\big\|^2
\,+\, \big\|\mcM_1(\bar{u}_1)\big\|^2.
\label{ptx1est9}
\end{equation}

Next, using that $ \; \big\|(\Dr -\alpha)v_1\big\|^2=\langle v_1\,,\,(\Dr -\alpha)^2v_1\rangle\;$ and $ \; (\Dr -\alpha)^2=(-\Dr)^2+2\al(-\Dr)+\al^2,\; $ we find
\begin{eqnarray}
\big\|(\Dr -\alpha)v_1\big\|^2
&=&
\|\Dr v_1\|^2
\; + \; 
2\, \al \,
\langle v_1 \,,\, (-\Dr) v_1\rangle
\; + \;
\al^2\,
\| v_1\|^2
\nonumber
\\[0.5em]
&=&
\big\|\Dr v_1\big\|^2\;+\;2\al\,\big\|\nabla^{\scriptscriptstyle rad}v_1\big\|^2
\; + \;
\al^2\,
\big\| v_1\big\|^2
\;.\label{Dmav1}
\end{eqnarray}
Using this relation and \eqref{e24225} in \eqref{ptx1est9} and using the inequality $ \|v_1\|^2_{H^1}\,\leq \,C_2\, Y_1, $ we obtain
\begin{equation}\label{pty1}
\p_tY_1 \;\leq\; 
-\, \frac14\; \big\|\Dr v_1\big\|^2
\;+\;
C_3\left(Y_1\,+\,Y_1^3 
\,+\, W\right)\,,
\end{equation}
where we collected all $L^2$-terms into $ C_3 Y_1, $ and where
\begin{equation}\label{y1rem}
W \equiv W(t) = \big\|u_1^\perp\big\|_{H^1}^4 \, + \,\big\|u_1^\perp\big\|_{H^1}^8\; .
\end{equation}
Eqs.~\eqref{pty1} and \eqref{y1rem} imply \eqref{vdifineq0} with $C=C_3$.
\end{proof} 

\begin{lemma}\label{lemdifin}
Let $ \; Y\equiv Y(t) \, $ be a positive and differentiable function satisfying the differential inequality
\begin{equation}\label{lemdifin1}
\partial_t Y \, \leq \,  C \,\left(Y + Y^3 + W\right) \,,
\end{equation}
with some constant $  C>0 , $ and a bounded function $ W\equiv W(t) $ and with initial condition $ \, Y(0) \geq 0$.  If $ \, Y(0)<1, $ then $ \, Y(t)\, $ satisfies the estimate
\begin{equation}\label{lemdifin2}
Y(t) \; \leq \; \left(Y(0)\,+\,\sup\limits_{0\leq s\leq t}W(s)\right)\;e^{C\,t} \;,
\end{equation}
for all $ \, t$'s such that $\, \big(Y(0)\,+\,\sup\limits_{0\leq s\leq t}W(s)\big)\,e^{C\,t} \leq 1.\,$
\end{lemma}

\begin{proof} Assume $ \, Y(t)\leq 1 \;,\forall t. $ We will justify this below. Then inequality \eqref{lemdifin1} implies
\begin{equation}\label{lemdifin3}
\partial_t Y \, \leq \, C \left(Y + W\right)\,.
\end{equation}Defining $\, Z(t) = e^{-C\,t}\,Y(t)\, $ and using \eqref{lemdifin3}, we obtain $ \, \p_t Z \leq e^{-C\,t}\,W.\, $ Integrating this gives $ \, Z(t)\,\leq\, Z(0) \,+\,C \int_0^t e^{-C\,s} \, W(s)\,ds\, ,\, $ which implies the `Duhamel inequality'
\begin{equation}\label{Duhin}
Y(t)\,\leq\, e^{C\,t}Y(0) \,+\, C \int_0^te^{C\,(t-s)}W(s)\,ds\,,
\end{equation}
which yields \eqref{lemdifin2}. 

For the time interval $ \, [0,T],\, $ where $\,T\,$ is such that $ \, \big(Y(0)\,+\,\sup\limits_{0\leq s\leq T}W(s)\big)\,e^{C\,T}  \,\leq \, 1,\, $ the assumption $ \, Y(t) \,\leq \, 1 \, $ holds. This proves the lemma.
\end{proof}

\begin{proof}[Proof of Theorem \ref{thm:average}]
\noindent We apply Lemma \ref{lemdifin} to the function $\,Y(t)=Y_1(t)  = - \tfrac12 \big\langle A_d\bm{v}(t) , \bm{v}(t)\big\rangle ,\, $ satisfying inequality \eqref{lemdifin1}, to obtain

\begin{equation}\label{X11dest}
Y_1(t) \;\leq\;\Big(Y_1(0)+ \sup\limits_{0\leq s\leq t}W(s)\Big)\;e^{C\,t}\;,
\end{equation}
where, by \eqref{vdifineq0} (Proposition \ref{vdifineq}),
\begin{equation}\label{Rsappl}
W(s):= 
\big\|u_1^\perp(s)\big\|_{H^1}^4
\;+\;
\big\|u_1^\perp(s)\big\|_{H^1}^8\;.
\end{equation}
By \eqref{rad-conv}, we have
\begin{equation}\label{Rest}
W(s) \,\leq\; 2 C^4  \, \big\|\bu_0- \bar\bu_{0} \big\|_{1,0}^4\; e^{-2\,\gamma\,\eps\,s}\;,
\end{equation}
which, together with \eqref{X11dest} yields
\begin{equation}\label{X11dest1}
Y_1(t)\;\leq\; \left(Y_1(0)+ 2C^4 \|\bu_0-\bar\bu_0 \|^4_{1,0}\right) \; e^{C\,t}\,,
\end{equation}
as long as $ \, \left(Y_1(0)+2C^4\|\bu_0-\bar\bu_0\|^4_{1,0}\right)\, e^{C\,t}<1. $ Finally, we recall that
$$ Y_1(t) \;\simeq\; \|\bm{v}(t)\|^2_{1,0} \;=\; \|\bar\bu(t)-\bm{w}(t)\|^2_{1,0} \;\,,
$$
which, together with \eqref{X11dest1}, implies \eqref{vest}, proving Theorem \ref{thm:average}.\end{proof}

\section{Concluding remarks} In this paper, we have shown that under some reasonable conditions, solutions of the FHN system on undulated cylindrical surfaces are exponentially approximated by their averages. The latter are close, for very long times, to solutions of the 1D system obtained from the standard FHN system by replacing $\partial_x^2$ by the radial surface Laplace-Beltrami operator. This possibly explains the effectiveness of the 1D (standard) FHN system. 

We expect that our techniques can be extended to warped cylindrical surfaces (with the graph function $\rho$ depending also on $\theta$). The results above make further steps along a long and challenging road to mathematical understanding of propagation of electrical impulses in neurons.

\paragraph*{\bfseries Acknowledgment} We are grateful to the anonymous referees for many useful remarks. The research work is implemented in the framework of H.F.R.I call “Basic research Financing (Horizontal support of all Sciences)” under the National Recovery and Resilience Plan “Greece 2.0” funded by the European Union – NextGenerationEU. (H.F.R.I. Project Number: 14910).
The research of I.M.S. is supported in part by NSERC Grant NA7901.

\vspace{0.3em}

\paragraph*{\bfseries Declarations}
\begin{itemize}

\vspace{-0.2em}

\item Conflict of interest: The authors have no conflicts of interest to declare that are relevant to the content of this article.
\item Data availability: Data sharing is not applicable to this article as no datasets were generated or analysed during the current study.
\end{itemize}

\appendix 
\section{Local existence}\label{sec:locglobexist}

Define the more general vector Sobolev spaces (cf. \eqref{eq:Hkl})
\begin{equation}\label{gvSs}
\vec{H}^{k,\ell}
:=\left\{\bu\in \vec L^2 :  
(-\Delta)^{k/2} u_1 \in L^2\,, \; \p_x^\ell u_2\in L^2 \right\}\;.
\end{equation}
We say, following \cite{TBS}, that Eq.~\eqref{FHNrho1} has a local \textit{mild} solution $  \,\bm{u}(t) \,$ if  $ \bm{u} \in C([0,T];H^{k,l}) $ for some $ T>0 $ (which depends on $\|\bm{u}_0\|_{k,l}$) and $ \bm{u}(t) $ satisfies the equation
\begin{equation}\label{inteqsol}
\bm{u}(t)=e^{tA}\bm{u}_0 + \int_0^t e^{(t-s)A} \mcN(\bm{u}(s))\,ds\;,
\end{equation}
where the operator $ A $ and the nonlinearity $ \mcN $ are defined in \eqref{A} and \eqref{N} and $ e^{tA} $ is the contraction semigroup on $ L^2 $ generated by $ A.$

Local well-posedness of Eq.~\eqref{FHNrho1} in $ \vec{H}^{2,1} $ in the sense of mild solutions is shown in \cite[Proposition 2.1]{TBS}. Arguing as in \cite{TBS}, we can prove the following:

\begin{proposition}\label{prop:localtime} For any $ \bu_0\in \vec{H}^{1,0} $ , system \eqref{FHNrho1} has a unique local mild solution in $ \, C([0,T],\vec{H}^{1,0})\, $ with
\begin{equation}\label{localtimelb}
T \;  \geq \; c \; \|\bu_0\|_{1,0}^{-2}
\end{equation}
for some constant $ c>0.$
\end{proposition}

\begin{proposition}\label{prop:strongsol}Under the conditions of Proposition \ref{prop:localtime}, the mild solution $ \, \bu(t)\, $ of \eqref{FHNrho1} is in fact a strong solution, i.e.
\begin{equation}\label{eq:stronglocsol}
\bu\in  C([0,T],\vec{H}^{1,0}) \cap C^1([0,T],\vec{H}^{-1,0}) \; .
\end{equation}
\end{proposition}

\begin{proof}[Proof of Proposition \ref{prop:strongsol}] We use that $ \, \bu(t) \, $ satisfies Eq.~\eqref{inteqsol},\, that $ \, e^{A t} \bu_0,\,\mcN(\bu(t))\in H^{1,0} \, $ and the integral in \eqref{inteqsol} is differentiable with respect to the upper limit, to obtain
\begin{equation}\label{eq:difinteq}
\partial_t \bm{u}(t)\,=\, A\, e^{A t}\bm{u}_0 \,+\, \int_0^t A \, e^{A (t-s)}\, \mcN(\bm{u}(s))\,ds\,+\,\mcN(\bu(t))\,.
\end{equation}
Now, since \;$  e^{At}\bm{u}_0 \, ,\, e^{A (t-s)} \mcN(\bm{u}(s)) \in H^{1,0} , \; $ we have that $ \; A e^{A t }\bm{u}_0 \,,\, A e^{A (t-s)} \mcN(\bm{u}(s)) \in H^{-1,0}\; $ and therefore, by \eqref{eq:difinteq}, $ \; \partial_t \bm{u}(t)\in H^{-1,0}.\, $ Similarly one shows the continuity of $ \, \partial_t\bu(t)\, $ in $ \, t.$
\end{proof}

\null\hfill\begin{tabular}[t]{l@{}}
\small 
Georgia Karali \\
\textit{E-mail address:} gkarali@math.uoa.gr
\\
\\
Konstantinos Tzirakis \\
\textit{E-mail address:} kostas.tzirakis@gmail.com
\\
\\
Israel Michael Sigal \\
\textit{E-mail address:} im.sigal@utoronto.ca
\end{tabular}
\end{document}